\documentclass[11pt]{article}
\usepackage{amssymb,amsfonts, amsmath,amsthm}
\usepackage[mathscr]{euscript}
\usepackage[utf8]{inputenc}
\usepackage{caption}
\usepackage{subcaption}
\usepackage{authblk}
\usepackage{color}
\usepackage{esvect}
\usepackage{empheq} 
\usepackage{float}
\usepackage{graphicx}
\usepackage{afterpage}
\usepackage{mwe}
\usepackage{comment}
\usepackage{mathtools}
\usepackage{xcolor}
\usepackage{epsfig,color,graphicx,graphics}
\usepackage{graphicx}

\numberwithin{equation}{section}
\allowdisplaybreaks
\newtheorem{thm}{Theorem}[section]
\newtheorem{defn}[thm]{Definition}
\newtheorem{cor}[thm]{Corollary}
\newtheorem{lemma}[thm]{Lemma}
\newtheorem{rmrk}[thm]{Remark}
\newtheorem{crl}[thm]{Corollary}
\newtheorem{pro}[thm]{Proposition}
\newtheorem{exmp}[thm]{Example}

\newcommand{\e}{\varepsilon}
\newcommand{\R}{\mathbb{R}}

\newcommand{\J}{\mathbb{J}}
\newtheorem{prop}[thm]{Proposition}

\newcommand{\abs}[1]{\left\vert{#1}\right\vert}

\newcommand{\ba}{\begin{array}}
\newcommand{\ea}{\end{array}}

\newcommand{\bthm}{\begin{thm}}
\newcommand{\ethm}{\end{thm}}
\newcommand{\bstp}{\begin{stp}}
\newcommand{\estp}{\end{stp}}
\newcommand{\blemma}{\begin{lemma}}
\newcommand{\elemma}{\end{lemma}}
\newcommand{\bprop}{\begin{prop}}
\newcommand{\eprop}{\end{prop}}
\newcommand{\bpf}{\begin{pf}}
\newcommand{\epf}{\end{pf}}
\newcommand{\bdefn}{\begin{defn}}

\newcommand{\edefn}{\end{defn}}
\newcommand{\brk}{\begin{rmrk}}
\newcommand{\erk}{\end{rmrk}}
\newcommand{\bcrl}{\begin{crl}}
\newcommand{\ecrl}{\end{crl}}
\newcommand{\beg}{\begin{exmp}}
\newcommand{\eeg}{\end{exmp}}
\newcommand{\norm}[1]{\left\|#1\right\|}

\newcommand{\beqn}{\begin{equation}}
\newcommand{\eeqn}{\end{equation}}

\newcommand{\pt}{\partial}

\renewcommand{\leq}{\leqslant}
\renewcommand{\geq}{\geqslant}

\newcommand{\beq}{\begin{equation}}
\newcommand{\eeq}{\end{equation}}
\newcommand{\bea}{\begin{eqnarray}}

\newcommand{\eea}{\end{eqnarray}}

\newcommand{\sh}{\mathcal{H}}

\newcommand{\loc}{\mathrm{loc}}

\newcommand{\indic}{\mathbf{1}}

\newcommand{\lb}[1]{{\color{purple!70!blue}#1}}

\def\({\left(}
\def\){\right)}
\def\I{{\mathcal I}}
\def\J{{\mathcal J}}
\def\lft{\text{left}}
\def\rght{\text{right}}

\def\rect{{\mathscr R}}
\def\ga{{\Gamma}}
\def\ba{{\bar\alpha}}
\def\bb{{\bar\beta}}

\title{}
\author

\begin{document}

\Large
 \begin{center}
\Large{Minimizing solutions of degenerate vector Allen-Cahn equations with three wells in $\R^2$}\\ 

\hspace{10pt}

\large
Lia Bronsard$^1$, \'Etienne Sandier$^2$, Peter Sternberg$^3$ \\

\hspace{10pt}

\small  
$^1$) McMaster University \\
bronsardmath@gmail.com
\\
$^2$) Universit\'e Paris-Est Cr\'eteil \\
sandier@u-pec.fr
\\
$^3$) Department of Mathematics, Indiana University, Bloomington \\
sternber@iu.edu

\end{center}

\hspace{10pt}

\normalsize

\begin{abstract}
We characterize all minimizers of the vector-valued Allen-Cahn equation in $\R^2$ under the assumption that the potential $W$ has three wells and that the associated  degenerate metric does not satisfy the usual strict triangle inequality, see \eqref{degenerate}. These minimizers depend on one variable only in a suitable coordinate system.

In particular, we show that no minimizing solutions to 
$ \Delta u=\nabla W(u)$ on $\R^2$
can approach the three distinct values of the potential wells.
\end{abstract}

\section{Introduction}

The study of entire solutions $u:\R^2\to \R^2$ to the vector Allen-Cahn system
\[
 \Delta u=\nabla W(u)
\]
 with a triple well potential $W:\R^2\to\R$ started with a first important result in \cite{BGS} and includes the more recent works \cite{Fusco,NAZG,ZG,EP} allowing for general, in particular, non-symmetric $W$. There is also a growing collection of work in the more general setting of $u:\R^n\to \R^m$ and multi-well  potentials $W:\R^m\to\R$, see \cite{AFS} and the references therein. A standard assumption in these works is that one has a strict triangle inequality between the potential wells for the degenerate metric $d(p,q)$ defined below in \eqref{Wdist}.
 Instead, here we pursue the degenerate case of equality as in \eqref{degenerate}, and under a set of further generic assumptions on $W$ presented below, we establish a rigidity result on the possible entire, minimizing solutions. 
 
 Let us state our assumptions on the potential. 
 
For $W\in C^3(\R^2;[0,\infty))$, we assume that
\[\{p\in \R^2: W(p)=0\}=P:=\{p_1,p_2,p_3\},\]
and we assume non-degeneracy of the potential wells in the sense that
\begin{equation}
\parbox{0.8\textwidth}{$D^2 W(p_\ell)$ has two {\em distinct, positive} eigenvalues for $\ell=1,2,3$.}\label{posdef}
\end{equation}

Additionally, we assume that for some $M>0$,
\begin{equation}
p\cdot \nabla W(p)\geq 0\quad\mbox{for}\;\abs{p}\geq M.
\label{Winfin}
\end{equation}

Next, we introduce the degenerate metric on $\R^2$ which arises in many studies of vector Allen-Cahn: 
\begin{align}
&d(p,q):=\nonumber\\
&\inf\left\{\sqrt{2}\int_{0}^{1}W^{1/2}(\gamma(t))
\abs{\gamma'(t)}\;dt:\gamma\in C^{1}([0,1],\R^{2}),
\;\gamma(0)=p,\;\gamma(1)=q\right\}.\label{Wdist}
\end{align}
We will assume that, contrary to the usual assumption of a strict triangle inequality, the following {\it equality} holds:
\begin{equation}
d(p_1,p_3) = d(p_1,p_2) + d(p_2,p_3).
\label{degenerate}
\end{equation}
This condition comes naturally in models such as tri-block co-polymers and leads to an interesting geometrical phenomenon, see \cite{ABLW}.
We will make the generic assumption that the geodesics from  $p_1$ to $p_2$ and from  $p_2$ to $p_3$ are {\it unique}. For the case when there is nonuniqueness of geodesics, see for example the work of \cite{ABG}. We also assume that the {\em only} length-minimizing geodesic between $p_1$ and $p_3$ goes through $p_2$.

We also note that an equivalent variational description of the distances $d(p_i,p_j)$ is given by
\beq
d(p_i,p_j)=\inf\left\{E(f,\R):\;f\in H^{1}_{loc}(\R,\R^{2}),
\;f(-\infty)=p_{i},\;f(\infty)=p_{j}\right\},\label{other}
\eeq
where 
\beq E(f,I) := \int_I \(\frac{1}{2}\abs{f'(t)}^2 + W(f(t))\)\;dt, \quad I\subset \R. \label{1denergy}\eeq
Then a minimizer $\zeta_{ij}$,  when it exists, is the parametrization of a geodesic joining $p_i$ and $p_j$. It satisfies the system  of ODE's
 \beq
\zeta_{ij}''(t)=\nabla W(\zeta_{ij}(t))\quad\mbox{for}\;-\infty<t<\infty,\;\zeta_{ij}(-\infty)=p_i,\;\zeta_{ij}(\infty)=p_j.\label{hetero}
\eeq 
From the perspective of ODE's, these parametrized geodesics $\zeta_{ij}$ represent heteroclinic connections between the potential wells. Our assumptions on $W$ imply that there exists a heteroclinic connection between $p_1$ and $p_2$ as well as between $p_2$ and $p_3$, but not between $p_1$ and $p_3$. We make the further assumptions that:  
 \beq\parbox{0.85\textwidth}{$\bullet$ The heteroclinic connections between $p_1$ and $p_2$ and between $p_2$ and $p_3$ are unique up to translation of the parameter\\ $\bullet$ The only nontrivial $H^1(\R)$ solution to the linearized equations $v''(t)= D^2W(\zeta_{ij})v$  is $v= \zeta_{ij}'$, for $ij = 12$ or $23$.)}
 \label{spectral}
 \eeq
This last hypothesis is formulated in a paper by M. Schatzman \cite{Sch} and proved there to be generic. We will crucially use Lemma~4.5 from \cite{Sch} in our proof of Lemma \ref{Sandier}. We will also need the following assumption, which is also certainly generic, although we will not try to prove this here:
 \beq\parbox{0.85\textwidth}{As $t\to +\infty$, we have $\zeta_{12}'/ |\zeta_{12}'|\to e$ and $\zeta_{32}'/ |\zeta_{32}'|\to -e$, where $e$ is an eigenvector of $D^2W(p_2)$ corresponding to the smallest eigenvalue.}
 \label{p2generic}
 \eeq

We remark that from the assumption \eqref{posdef}, it follows that each $\zeta_{ij}$ approaches its end-states $p_i$ and $p_j$ at an exponential rate, i.e.
\beq
\abs{\zeta_{ij}(s)-p_j}\leq Ce^{-cs}\quad\mbox{as} \;s\to\infty\label{expdecay}
\eeq
for some constant $c$ depending on $W$,  with a similar estimate holding as $s\to -\infty$. 
(See Section 2 for details.)

Next,  for any bounded open $\Omega\subset\R^2$ and any $u\in H^1(\Omega,\R^2)$, we define
\begin{equation}
    E(u,\Omega)=\int_{\Omega}\(\frac{1}{2}\abs{\nabla u}^{2}+ W(u)\)\;dx\,dy.
    \label{energy}
\end{equation}
We say that $u:\R^2\to\R^2$ is a minimizer of $E$ if it minimizes $E$ with respect to its own boundary conditions on any bounded domain $\Omega\subset \R^2$. We remark that such an entire solution is alternatively referred to as a local minimizer (in the sense of De Giorgi) in much of the literature.

Our main result is the following:

 \begin{thm}\label{main}
 Assume $u:\R^2\to\R^2$ is a nonconstant minimizer of $E$. Then, possibly after rotating and translating the coordinates, either
$u(x,y) \equiv \zeta_{12}(x)$ or 
    $u(x,y) \equiv \zeta_{23}(x)$. 
 \end{thm}
 
Our article is organized as follows: in Section 2, we present some preliminary results, in particular we show that, up to a translation and/or rotation, the only possible blowdown limits of a nonconstant minimizer $u$ are the piecewise constant functions $H_{12}, H_{23}$ and $H_{13}$,  where
\begin{equation}
    H_{ij}(x,y) := p_i\indic_{\{x<0\}}(x,y) + p_j\indic_{\{x>0\}}(x,y).
    \label{half}
\end{equation}
 A result of \cite{EP} shows that if, up to translation and rotation, the blowdown limits of $u$ coincide with those of either $\zeta_{12}$ or $\zeta_{23}$, namely $H_{12}$ or $H_{23}$ respectively, then $u$ itself must be equal to $\zeta_{12}$ or $\zeta_{23}$ (up to the invariances of the problem). The proof then consists of ruling out the possibility that  the blowdown limit is $H_{13}.$ 

The crucial fact we use is that under our hypothesis on $W$, for $u:\R\to\R^2$ to transition optimally from $p_1$ to $p_3$, there must be two  transitions: one from $p_1$ to $p_2$ and one from $p_2$ to $p_3$, but these two transitions must also  be as far apart as possible. We make this statement precise in Section 3, where we study the structure of almost minimizers in $1D$. In Section 4, we obtain a precise upper bound for the energy of minimizers on the ball $B_R$ for $R\gg 1$, which exhibit nearly optimal $p_1$-$p_3$ transitions on the boundary. Finally, in Section 5, we prove a lower bound for the energy on $B_R$ of a minimizer whose blowdown is $H_{13}$, and observe that this lower bound is incompatible with the upper bound we computed before, and hence conclude the proof of the theorem.

\section*{Acknowledgements}
The authors would like to thank A.Monteil for his help in proving Lemmas~\ref{lem20} and \ref{lem21}.

This research benefited from invitations of L.B and P.S to Université Paris-Est Créteil and of E.S. to McMaster University, that they wish to thank for their hospitality. These visits were made possible by funding from an International Emerging Action of INSMI-CNRS of E.S. and the NSERC discovery grant of L.B.
The research of P.S. was supported by a Simons Collaboration grant 585520 and an NSF grant DMS 2106516.

\section{Preliminaries and background}

We begin with some preliminary results on the behavior of a geodesic $\gamma$ connecting a point $a$ to a well $p\in P$. Though, of course, a distance-minimizing geodesic in the plane satisfies a second order system of ODE's, one can alternatively express its dynamics as a gradient flow with respect to the distance $d$ given by \eqref{Wdist}. Thus, we have
\begin{equation}\gamma(0) = a,\quad \gamma(+\infty) = p,\quad\text{and}\quad \gamma'(t) =-\nabla d_p(\gamma(t)),\label{linzz}
\end{equation}
where $d_p(\cdot) = d(p,\cdot)$. In the case where $W$ is a quadratic in a neighborhood of $p$, that is, where \[W(q)\equiv \lambda_1(q^{(1)}-p^{(1)})^2+ \lambda_2(q^{(2)}-p^{(2)})^2
\] 
in an orthonormal coordinate system given by the eigenvectors of $D^2W(p)$,
it is shown in \cite{ABCDS}, Section 2.1, that $d_p$ is given by 
\[
d_p(q)=\frac{1}{2}\left(\sqrt{\lambda_1}(q^{(1)}-p^{(1)})^2+\sqrt{\lambda_2}(q^{(2)}-p^{(2)})^2\right),
\]
so that in the purely quadratic case, with say $p=0$, one has 
\begin{equation}
\left(\gamma^{(1)},\gamma^{(2)}\right)'=-\left(\sqrt{\lambda_1}\gamma^{(1)},\sqrt{\lambda_2}\gamma^{(2)}\right).\label{linz}
\end{equation}


Therefore if $p$ is rest point of the ODE \eqref{linzz}, then the linearization of the ODE at $p$ is necessarily \eqref{linz}, and consequently, it is a contraction. It then follows from \cite{hartman} that for some $\eta>0$, there exists a $C^1$ diffeomorphism $h_p$ from $B_\eta(p)$ to a neighborhood of $0\in\R^2$ such that  if $\gamma' =-\nabla d_p(\gamma)$, with $\gamma$ taking values in $B_\eta(p)$, then $X = h\circ\gamma$ solves $X' = AX$ with $A:=\,-{\rm{diag}} [\sqrt{\lambda_1},\sqrt{\lambda_2}]$. 

Note that this result does not hold in $\R^n$ for $n>2$ (see \cite{hartman}), however in that case there still exists a homeomorphism which is differentiable at $p$, while its inverse is differentiable at $0$. This result is proved in \cite{zlz} for discrete dynamical systems but, as mentionned in \cite{hartman}, using  \cite{shlomo}, Lemma~4, this carries over to smooth dynamical systems. 

We thus have
\begin{lemma}\label{local} Assume that $W:\R^2\to\R_+$ is $C^2$ and that $\nabla W(p) = 0$, $D^2W(p) > 0$. Then there exists $\eta>0$ and $C^1$-diffemorphism $h$ from $U = B_\eta(p)$ to a neighbourhood of $0\in\R^2$ such that $Dh(p) = \text{Id}$ and such that if $\gamma:\R_+\to U$ is a geodesic with the property that  $\gamma(t)\to p$ as $t\to +\infty$,  then $h\circ u = X$, where $X$ solves $X' =\,-{\rm{diag}} [\sqrt{\lambda_1},\sqrt{\lambda_2}]X$.

    In particular, choosing  orthonormal coordinates $(x,y)$ in which $p=0$ and $\frac12 D^2W(p) = \begin{pmatrix}
    \lambda_1 & 0\\
    0 & \lambda_2\end{pmatrix}$, and assuming  $0<\lambda_1 < \lambda_2$,  it follows that if $\gamma(0) = (x_0,y_0)$, then $\gamma(t)\approx (x_0 e^{-\sqrt{\lambda_1} t}, y_0 e^{-\sqrt{\lambda_2} t})$ as $t\to +\infty$.
\end{lemma}

The property \eqref{expdecay}  also follows from this lemma. 
  
Next, we note that even though minimizing \eqref{1denergy} fixes a parametrization of the geodesic, there is still translation invariance. We may remove this invariance and then unambiguously define $\zeta_{ij}:\R\to\R^2$ by requiring that
 \beq 
 d(p_i,\zeta_{ij}(0)) = d(p_j,\zeta_{ij}(0)).
 \label{middle}
 \eeq
From now on we will write
\beq 
 d_{ij} = d(p_i,p_j).
 \label{dij}
 \eeq

We now summarize the well-known Gamma-convergence properties of the rescaled energy $E_R(u,\Omega)$ on a bounded smooth domain $\Omega$ where
\begin{equation}
E_R(u,\Omega):=\int_{\Omega}\(\frac{1}{2R}\abs{\nabla u}^{2}+ RW(u)\)
\;dx\,dy.\label{Renergy}
\end{equation}

Then the $L^1$  Gamma-limit of the sequence $\{E_R\}$ is the functional $E_0$ given by
\beq
E_0(u,\Omega) := \sum_{1\leq i<j\leq 3}d_{ij}\sh^{1}(\partial^*S_{i}\cap\partial^*S_{j} \cap \Omega),\label{ez}
\eeq
defined on $BV_\loc(\R^2,P)$, where $P = \{p_1,p_2,p_3\}$, where $S_{j}:=u^{-1}(p_{j})$ for $j=1,2,3$, and $\partial^*S$ refers to the reduced boundary of a set
$S$ of finite perimeter, cf. \cite{giusti}. This follows from \cite{baldo} in the absence of boundary conditions, and from \cite{Gaz} for the case of Dirichlet boundary conditions. Note that if $u$ is a minimizer of $E_0$, and in the case were the triangle inequality is strict for the distance between the wells, the interface between the sets $S_j$ consists of line segments that can meet at triple junctions, the angles they make  are determined through stationarity by the $d_{ij}$'s and are nonzero. 

Next, for any $R>0$, and any $u:\R^2\to\R^2$, we introduce the notation $u_R=u(R\,\cdot)$ to denote the {\it blowdown} of the function $u$. One can then readily check that the relation $E(u,R\Omega)=RE_R(u_R,\Omega)$ holds. We will make use of the following well known fact in the proof of the proposition to follow.

 \begin{rmrk}\label{locunif}
    If $u:\R^2\to\R^2$ is a minimizer of $E$, and if  $u_R\to u_0$ in $L^1_\loc$, where $u_0\in BV_\loc(\R^2,P)$, then $u_R$ converges locally uniformly to $p_j$ in the interior of each $S_{j}:=u_0^{-1}(p_{j})$. See e.g. Prop. 4.2 of \cite{EP}.
 \end{rmrk}

Recalling the definition of $H_{ij}$ given in \eqref{half},  we have  the following characterization of the possible limit of blowdowns.

\begin{pro} \label{gamma}Assume $u:\R^2\to\R^2$ is a non-constant minimizer of $E$, then for any sequence $\{R_n\}$  of radii tending to $+\infty$, there exists a subsequence, still denoted $\{R_n\}$, such that, possibly after rotating the coordinates, $u(R_n\;\cdot)$ converges to  $H_{ij}$ in $L^1_\loc$, for some $1\le i\neq j\le 3$. 

    Moreover, for any bounded open $\Omega\subset \R^2$, it holds that 
    $$E_{R_n}(u(R_n\;\cdot),\Omega) = E_0(H_{ij},\Omega)+o(1).$$
\end{pro}

\begin{proof} Applying Proposition 3.1 of \cite{EP}, which in particular does not rely on a strict triangle inequality, we have that for any sequence $\{R_n\}\to\infty$, there exists a subsequence, still denoted by $\{R_{n}\}$, and a function $u_0\in BV_{{\rm loc}}\left(\R^2;P\right)$ such that the blowdowns $\{u_{R_{n}} \}$ of $u$ satisfy
\beq
u_{R_{n}}\to u_0\quad\mbox{in}\; L^1_{{\rm loc}}(\R^2;\R^2).\label{Lone}
\eeq
Furthermore, $u_0$ is a minimizer of $E_0$ and in light of the Gamma-convergence of $E_R$ to $E_0$, we have that 
 \begin{equation}
 E_{R_n}(u(R_n\;\cdot),\Omega) = E_0(u_0,\Omega)+o(1)\label{quinn}
 \end{equation}
 for any bounded open $\Omega\subset \R^2$. 
 
 We first eliminate the possibility that $u_0=p_j$ for some $j\in\{1,2,3\}.$ Suppose, by way of contradiction, that say $u_0=p_1.$ Then by \eqref{quinn}, we have that $E_{R_n}(u_{R_n},B_1)=o(1)$ as $n\to\infty$, or equivalently, $E(u,B_{R_n})=o(R).$ Hence, an application of the Mean Value Theorem leads to a sequence of radii, say $\rho_n\to\infty$, such that
 $E(u,\partial B_{\rho_n})\to 0.$ Also, through an appeal to Proposition 4.2 of \cite{EP}, we have that
 \[
 u_{\rho_n}\to p_1\quad\text{locally uniformly on}\;\R^2,
 \]
 which implies, in particular, that
 \begin{equation}
 \sup_{(x,y)\in \partial B_{\rho_n}}|u(x,y)-p_1|=o(1).\label{rosalyn}
 \end{equation}
With \eqref{quinn} and \eqref{rosalyn} in hand, we may construct a low-energy competitor on the ball $B_{\rho_n}$ by linearly interpolating between the values of $u$ on $\partial B_{\rho_n}$ and $p_1$ on $\partial B_{\rho_n-1}$. A straightforward estimation reveals that the cost on this annulus is $o(1)$. (One can view the estimates in \cite{EP}, pages 4179-4181 for a similar estimation in a more complicated situation.) Then filling in the ball $ B_{\rho_n-1}$ with the value $p_1$ we have shown that the minimizer $u$ satisfies $E(u, B_{\rho_n})=o(1)$. Hence, in fact  $E(u,\R^2)=0$ and necessarily $u\equiv p_1$, contradicting the assumption that $u$ is nonconstant.

Next we appeal to the regularity result for minimizing partitions of the plane. Appealing to the result \cite{White}, Thm. 3, it is argued in \cite{Morgan}, Thm. 4.3 that even in the absence of a strict triangle inequality for, say, the cost $d(p_1,p_3)$, the phase boundaries of any minimizing partition of the plane consist of a union of line segments joined at triple junctions. What is more, in light of \eqref{degenerate}, there can be no junctions at all since stationarity at such a junction would require the phase corresponding to $p_3$ to be enclosed with an angle of size $0$.  It follows that the only possibility for a minimizing partition would be a collection of parallel strips--that is, phase boundaries consisting of parallel lines. 

We may assume without loss of generality that these lines are vertical, so that $u_0(x,y) = u_0(x)$. It is then not difficult to argue by constructing competitors that $u_0$, as a function of one variable, must be a minimizer for the one-dimensional transition energy. Then a simple examination of this problem shows that if $u_0$ has more than one transition, then the only possibility is that it transitions from $p_1$ to $p_2$ and then from $p_2$ to $p_3$, either from left-to-right or from right-to-left. 

Also, going back to the 2-dimensional picture, if $u_0$ is not a constant, then one of the transition lines must go through the origin, for otherwise we would have $u_0\equiv p_i$ in $B(0,\eta)$, which would imply, from the uniform convergence of $u_{R_n}$ to $u_0$ on $B(0,\eta)$, that $u\equiv p_i$. 

Thus, with no loss of generality, suppose the phase boundaries are vertical and that the $y$-axis separates the $p_1$-phase from the $p_2$-phase. Then there are numbers $a<0$ and $b>0$ such that $u_0=p_1$ for $a<x<0$ and $u_0=p_2$ for $0<x<b$ with $i\not= k.$ Let $\{\eta_j\}\to 0$ be an arbitrary sequence of positive numbers. Then for each $j$ there exists a value $R_{n_j}$ (which we denote simply by $R_j$) such that $\eta_j^{1/3}R_j\to\infty$ and such that
\begin{equation}
\norm{u_{R_j}-p_1}_{L^1\big([a/2,-\eta_j]\times [-1/\eta_j,1/\eta_j]\big)}<\eta_j
\end{equation}
and
\begin{equation}
\norm{u_{R_j}-p_2}_{L^1\big([\eta_j,b/2]\times [ -1/\eta_j,1/\eta_j]\big)}<\eta_j
\end{equation}
But letting $\tilde{x}=x/\eta_j^{1/3}$ and $\tilde{y}=y/\eta_j^{1/3}$ these two conditions are equivalent to the conditions
\begin{equation}
\norm{u_{\eta_j^{1/3}R_j}-p_1}_{L^1\big([a/2\eta_j^{1/3},-\eta_j^{2/3}]\times [-1/\eta_j^{4/3},1/\eta_j^{4/3}]\big)}<\eta_j^{1/3}
\end{equation}
and
\begin{equation}
\norm{u_{\eta_j^{1/3}R_j}-p_2}_{L^1\big([\eta_j^{2/3},b/2\eta_j^{1/3}]\times [-1/\eta_j^{4/3},1/\eta_j^{4/3}]\big)}<\eta_j^{1/3}.
\end{equation}
In other words, along the sequence $\{\eta_j^{1/3}R_j\}\to\infty$ the blowdowns converge in $L^1_{loc}$ to $H_{12}$. Consequently, we may apply Theorem 3.1 of \cite{EP} to conclude that $u\equiv \zeta_{12}$ up to a translation. But this contradicts the assumption that the blowdown consists of more than one phase boundary. Of course, the same argument would apply if phases $p_2$ and $p_3$ were adjacent.
\end{proof}

 Ultimately, we will show that the blowdown limits can only be $H_{12}$ or $H_{23}$, not $H_{13}$ and then that the only minimizers are $\zeta_{12}$ and $\zeta_{23}$, up to the invariances of $E$, thus proving Theorem \ref{main}.


\section{Structure of almost minimizers in 1D}

\begin{lemma}
\label{lem20}
Assume $0$ is a nondegenerate minimizer of $W:\R^2\to\R_+$. Then there exists $\eta>0$ and $R_0,C, c>0$ such that if  $R>R_0$, the following holds. For any $u:[0,R]\to\R^2$ such that $u(0) = a$, with $a$ such that $d(0,a)\le \eta$ it holds that $$E(u)\ge d(0,a) - Ce^{-cR}.$$
\end{lemma}

The proof uses a strategy suggested to us by A.Monteil (\cite{antonin}), which is used again in the next lemma.

\begin{proof}
First we note that we may assume that $u$ minimizes the energy on $[0,R]$ with respect to the initial condition $u(0)=a$. In particular, choosing $\eta>0$ sufficiently small so that $W$ is convex on the set of points \lb{ $p$} satisfying $d(0,p)\le \eta$, we may assume that $d(u(t),a)\le \eta$ for any $t$. Let $\gamma_a:[0,+\infty)\to \R^2$  be the distance minimizing geodesic from $a$ to $0$, parametrized so that $E(\gamma_a) = d(a,0)$. Let also $p = u(R)$ and $q = \gamma_a(R)$.

We define $\tilde u:\R_+\to\R^2$ by $\tilde u(t) = u(t)$ for $t\in[0,R]$ and $\tilde u(t) =\gamma_{p}(t-R)$ for $t\ge R$, where $\gamma_p:[0,+\infty)\to \R^2$  is the distance minimizing geodesic from $p$ to $0$ such that $E(\gamma_p) = d(p,0).$

The maps $\tilde u$ and ${\gamma}_a$ have the same boundary conditions on $[0,+\infty]$, therefore, letting $Q$ be the quadratic form associated to $\frac{1}{2}D^2W(0)$, we have 
\begin{equation}
\begin{split}
E(\tilde u) - E(\gamma_a) &=\int_{0}^{+\infty}\(\frac12 |\tilde u'|^2 + W(\tilde u)\) - \(\frac12 |\gamma_a'|^2 + W(\gamma_a)\) \;dt\\
&\ge \int_{0}^{+\infty}\frac12\( |\tilde u'|^2 - |\gamma_a'|^2\) + DW(\gamma_a)(\tilde u-\gamma_a) + (1-C\eta) Q(\tilde u-\gamma_a) \;dt\\
&= \int_{0}^{+\infty}\gamma_a'\cdot (\tilde u-\gamma_a)'+ \frac12 |(\tilde u-\gamma_a)'|^2 +
DW(\gamma_a)(\tilde u-\gamma_a)\\& + (1-C\eta) Q(\tilde u-\gamma_a) \;dt\\
&= \int_{0}^{+\infty}\frac12 |(\tilde u-\gamma_a)'|^2 + (1-C\eta) Q(\tilde u-\gamma_a) \;dt,
\end{split}
\end{equation}
where we have used the fact that $\gamma_a'' = DW(\gamma_a)$.
Therefore (and recalling that $C$ denote a constant that may change from line to line)
$$E(\tilde u) - d(a,0)\ge \sqrt{1-C\eta}\,\sqrt{2}\int_{0}^{+\infty}\sqrt{Q(\tilde u-\gamma_a)} |(\tilde u-\gamma_a)'| \;dt.$$

We now introduce the distance associated to the form $Q$, which we denote by $d_Q,$ defined by replacing $W$ by $Q$ in \eqref{Wdist}. Then  the inequality above can be phrased as 
\begin{equation*}
E(\tilde u) - d(a,0)\ge 2\sqrt{1-C\eta}\,d_Q(0,p-q),
\end{equation*}
where the factor of $2$ emerges since the path $\tilde u-\gamma_a$ goes from $0$ to $u(R)-\gamma_a(R) = p-q$ and then returns. Since $E(\tilde u) = E(u)+d(p,0)$ and $d(a,0) = d(a,q) + d(q,0)$ we deduce that
\begin{equation}
E(u) \ge d(a,q) + d(q, 0) -d(p,0) + 2\sqrt{1-C\eta}\,d_Q(0,p-q).\label{thing}
\end{equation}
Now (see \cite{ABCDS}, Remark 2.4) we have that $x\to \sqrt{d_Q(x,0)}$ is a norm on $\R^2$, that we denote $\|x\|_Q$, therefore 
$$d_Q(0,p-q) = \|p-q\|_Q^2 \ge d_Q(0,p) + d_Q(0,q) - 2\|p\|_Q\|q\|_Q \ge d_Q(p,q) - 2\|p\|_Q\|q\|_Q.$$
Together with \eqref{thing}, and using the fact that $d_Q\ge (1-C\eta)d$ this implies that
$$E(u)\ge d(a,q) + d(q, 0) -d(p,0) + d(p,q) + (1-C\eta) \|p-q\|_Q^2 - 4\|p\|_Q\|q\|_Q.$$
and therefore,
\begin{equation}
E(u)\ge d(a,q) + (1-C\eta) \|p-q\|_Q^2 - 4\|p\|_Q\|q\|_Q.\label{ugh}
\end{equation}
Now consider $\eta$  chosen sufficiently small such that $C\eta<1/2$. Then we claim that \eqref{ugh} implies the inequality
\[
E(u)\ge d(a,q) - C\|q\|_Q^2.
\]
Indeed if, say  $\|p\|_Q>M\|q\|_Q$ for some large enough $M$, then in fact \eqref{ugh} would imply that $E(u)\ge d(a,q)$ since the middle term on the right-hand side of \eqref{ugh} would then dominate the last term. On the other hand, if $\|p\|_Q\leq M\|q\|_Q$,
then one has the claim with $C=4M$ by ignoring the middle term on the right-hand side of \eqref{ugh}.

Thus, 
\[
E(u)\ge d(a,q)-Cd(0,q)=d(0,a)-d(0,q)-C\|q\|_Q^2=d(0,a)-Cd(0,q).
\]
Since $\gamma_a(t)$ converges to $0$ exponentially fast as $t\to +\infty$, we find that $d(0,q) = d(0,\gamma_a(R))\le Ce^{-cR}$, and so
$$E(u) \ge d(0,a) - Ce^{-cR}.$$
\end{proof}

\begin{lemma}\label{lem21} 
 Assume $0$ is a nondegenerate minimizer of $W:\R^2\to\R_+$ and choose coordinates $(x,y)$ such that $D^2W(0) = \begin{pmatrix}
    \lambda_1 & 0\\
    0 & \lambda_2\end{pmatrix}$ with $0<\lambda_1 < \lambda_2$ in these coordinates. Let $d$ be the distance associated to $W$ as in \eqref{Wdist}.     For  $\eta>0$ and  $\e>0$ sufficiently small depending on $W$, if $f\in H^1([-R,R],\R^2)$  satisfies $f(-R)=a$ and $f(R)=b$ with $a=(a_1,a_2)$ and $b=(b_1,b_2)$, where $d(a,0)=\eta=d(b,0)$  
    $a_1<0$, $b_1>0$, and $|a_2|,|b_2|\le \e |a_1|, \e |b_1|$,
    then one has 
    \begin{equation}
        E(f,[-R,R])\ge 2\eta +c\left(\min_{t \in[-R,R]}|f(t)-p_2|^2+ e^{-CR}\right),
    \end{equation}
    for some $c,C>0$ independent of $R$.
\end{lemma}
\begin{proof}
Let $f\in H^1([-R,R],\R^2)$ be a minimizer of $E$ from \eqref{1denergy} with $I=[-R,R]$ under the boundary conditions $f(-R)=a$ and $f(R)=b$. Then  for $\eta$ sufficiently small depending on $W$, we have $d(f(t),0)\le \eta$ for any $t\in(-R,R)$, and furthermore  there exists $t_0\in (-R,R)$ such that $f(t_0) = q := (0,q_2)$, since $a_1b_1<0$. 

We denote by $\gamma_p:[0,+\infty)\to \R^2$ the distance minimizing geodesic from any point $p$ to $0$, parametrized so that $E(\gamma_p,[0,+\infty)) = d(p,0)$. We define $f_a:[-R,+\infty)\to\R^2$ by requiring that 
$$f_a(t) = \begin{cases} f(t) &\mbox{if}\;t\in[-R,t_0],\\
\gamma_{q}(t-t_0) &\mbox{if}\;t\in[t_0,+\infty).\end{cases}$$
The maps $f_a$ and $\tilde{\gamma}_a:=\gamma_a(\cdot +R)$ have the same boundary conditions on $[-R,+\infty]$, therefore 
\begin{equation}
\begin{split}
E(f_a,[-R,\infty)) - & E(\gamma_a(\cdot+R),[-R,\infty)) \\ &=\int_{-R}^{+\infty}\(\frac12 |f_a'|^2 + W(f_a)\) - \(\frac12 |\tilde{\gamma}_a'|^2 + W(\tilde\gamma_a)\) \;dt\\
&\ge \int_{-R}^{+\infty}\frac12\( |f_a'|^2 - |\tilde\gamma_a'|^2\) + DW(\tilde\gamma_a)(f_a-\tilde\gamma_a) + (1-C\eta) Q(f_a-\tilde\gamma_a) \;dt\\
&= \int_{-R}^{+\infty}\tilde\gamma_a'\cdot (f_a-\tilde\gamma_a)'+ \frac12 |(f_a-\tilde\gamma_a)'|^2 +
DW(\tilde\gamma_a)(f_a-\tilde\gamma_a) \\ &+ (1-C\eta) Q(f_a-\tilde\gamma_a) \;dt\\
&= \int_{-R}^{+\infty}\frac12 |(f_a-\tilde\gamma_a)'|^2 + (1-C\eta) Q(f_a-\tilde\gamma_a) \;dt,
\end{split}
\end{equation}
where $Q$ is the quadratic form associated to $\frac{1}{2}D^2W(0)$, and where we have used the fact that $\tilde\gamma_a'' = DW(\tilde\gamma_a)$.

Therefore 
$$E(f_a,[-R,\infty)) - d(a,0)\ge \sqrt{2}(1-C\eta)\int_{-R}^{+\infty}\sqrt{Q(f_a-\tilde\gamma_a)} |(f_a-\tilde\gamma_a)'| \;dt.$$

We now recall the distance associated to the form $Q$, which we denote by $d_Q,$ analogously to \eqref{Wdist}. Then  the inequality above can be phrased as 
\begin{equation}
E(f_a) - d(a,0)\ge 2(1-C\eta)d_Q(0,q-\tilde\gamma_a(t_0)),\label{lia}
\end{equation}
where the factor of $2$ emerges since the path $f_a-\tilde\gamma_a$ goes from $0$ to $q-\tilde\gamma_a(t_0)$ and then returns. Denoting the coordinates of $\tilde\gamma_a(t_0)$ by $(\alpha_1,\alpha_2)$, one has that 
\begin{equation}
d_Q(0,q-\tilde\gamma_a(t_0))=\frac{1}{\sqrt{2}}\big(\sqrt{\lambda_1}\alpha_1^2+\sqrt{\lambda_2}(q_2-\alpha_2)^2\big)\quad\text{while}\quad d_Q(0,q)=\frac{1}{\sqrt{2}}\sqrt{\lambda_2}\,q_2^2\label{etienne}
\end{equation}
(again, see \cite{ABCDS}, Remark 2.4).
Moreover, from Lemma~\ref{local} and in view of the hypothesis $|a_2|\le \e |a_1|,$ it holds that $|\alpha_2|\le C\e|\alpha_1|$. Since $0<\lambda_1 <\lambda_2$ and since $d(0,q)\le (1+C\eta)d_Q(0,q)$, one can then use the inequalities \eqref{lia} and \eqref{etienne} to find that for $\eta$ and $\e$ chosen small enough
\[
2(1-C\eta)d_Q(0,q-\tilde\gamma_a(t_0))-d(0,q)\ge c\(\alpha_1^2+q_2^2\)\ge c\(d(0,\tilde\gamma_a(t_0))+|f(t_0)-p|^2\),
\]
by considering separately the case where $|q_2|$ is negligible compared with $|\alpha_2|$ and where it is not, and using that $|\alpha_2|\le C \e|\alpha_1|$ in the second inequality. Here $c$ is a positive constant depending on $W$.

Inserting this inequality into \eqref{lia} we obtain
\begin{equation}
E(f_a)=E(f,[-R,t_0])+d(0,q)\ge d(0,a)+d(0,q)+c\(d(0,\tilde\gamma_a(t_0))+|f(t_0)-0|^2\).
\end{equation}
Since $t_0+R\le 2R$, we have, using Lemma~\ref{local}, that  $d(0,\tilde\gamma_a(t_0))\ge ce^{-CR}$ and therefore,
\[
E(f,[-R,t_0])\ge d(0,a)+c\(|f(t_0)-p|^2+e^{-CR}\).
\]
By a similar argument we have the inequality
\[
E(f,[t_0,R])\ge d(0,b)+c\(|f(t_0)-p|^2+e^{-CR}\).
\]
Adding these two inequalities proves the lemma.

\end{proof}

\begin{lemma} \label{Sandier}There exists $\eta>0$ and $C>0$ such that for any $\gamma>0$ small enough and any $R>0$ large enough the following holds.

Assume $f:[-R,R]\to\R^2$ is such that 
\begin{itemize}
\item for every $t\in[-R,-R/2]$, $|f(t)-p_1|\le \eta$,
\item for every $t\in[R/2,R]$, $|f(t)-p_3|\le \eta$,
\item $E(f,[-R,R])\le d_{13}+\gamma$.
\end{itemize}
Then, there exist $a$, $b$  and $T$, with $-R/2<a<T<b< R/2$  such that 
\begin{itemize}

    \item $\|f(a+\cdot) - \zeta_{12}(\cdot)\|_{H^1([-R/2,T])}^2\le C\big(\gamma+e^{-cR}\big)$,
    \item $\|f(b+\cdot) - \zeta_{23}(\cdot)\|_{H^1([T,R/2])}^2\le C\big(\gamma+e^{-cR}\big)$,
    \item $|f(T)-p_2|^2< C\big(\gamma+e^{-cR}\big)$
    \end{itemize} 
\end{lemma}

\begin{proof} Let $\eta>0$ be such that $W$ is strictly convex in $B(p_i,3\eta)$ for $i=1,2,3$. 
    It suffices to show that if $R_n$ converges to $+\infty$ and $\gamma_n$ converges to $0$, and if $f_n$ satisfies the above hypotheses, then for any $n$ large enough there exist $a_n$, $b_n$ and $T_n$ with $-R_n/2<a_n<T_n<b_n<R_n/2$
    such that the conclusion of the proposition holds for some $C>0$ independent of $n$. Throughout the proof, $C$ denotes a positive large constant independent of $n$, and $c$ denotes a positive small constant independent of $n$.

    Let $A_n$ be the set of $t\in[-R_n,R_n]$ such that $d(f_n(t),P)>\eta$. Then $A_n$ is a union of disjoint open intervals included in $(- R_n/2,R_n/2)$. We denote $\I = \{I_1,\dots,I_l\}$ the family of those intervals in which there exists $t$ such that $d(f_n(t),P)>4\eta$, and call them the transition intervals. If $I$ is a transition interval, then $E(f_n,I)\ge 3\eta$, using the trivial lower bound $E(f_n,[t,t'])\ge d(f_n(t),f_n(t'))$ which results from the definition of $d$.

    Because of the energy bound, the number of transition intervals is thus bounded independently of $n$. It is easy to check that the energy bound also implies that their size is bounded above and below independently of $n$. 
    
    The complement of $\cup_j I_j$ is a union of the sequence of closed disjoint intervals $\J = \{J_0,\dots,J_k\}$. On each $J_i$, with $0\le i\le k$, the map $f_n$ is close to one of the wells that we denote $q_i$. For $0\le i\le k$ we let $J_i = (x_i,y_i)$, so that $x_0 = -R_n$ and $y_k = R_n$, and so that  $d(f_n(x_i), q_{i}) \le \eta$ and $d(f_n(y_i), q_{i}) \le \eta$.
    we say that $I_i$ is a transition interval from the well $q_{i-1}$ to the well $q_i$. Note that $q_0 = p_1$ and $q_k = p_3$.

    From Lemma~\ref{lem20}, we have $E(f_n,J_0)\ge d(p_1,f_n(y_1)) - c e^{-CR_n}$ and $E(f_n,J_k)\ge d(f_n(x_k),p_3) - c e^{-CR_n}.$

    We claim that, in the sequence $q_0,\dots,q_k$, no point in $P$ can appear twice. Indeed, if we had $q_i = q_j = p\in P$ with $i<j$, then we would have $d(f_n(y_i),p)\le\eta$ and $d(f_n(x_j),p)\le\eta$. By using the trivial lower bound $E(f,[t,t'])\ge d(f(t),f(t'))$ on each interval in $\I$ or $\J$ to the left of $y_i$, including $J_i$ we find, denoting $I_\lft$ their union, 
    $$E(f_n,I_\lft )\ge d(p_1,f_n(y_i)) - c e^{-CR_n}.$$
    Similarly, denoting $I_\rght$ the union of the intervals to the right of $x_j$, including $J_j$, we have
    $$E(f_n,I_\rght )\ge d(f_n(x_j),p_3) - c e^{-CR_n}.$$
    But between $y_i$ and $x_j$, there is at least a transition interval, on which the energy is bounded below (as we have seen) by $3\eta$. Thus 
    $$E(f_n, [-R_n,R_n])\ge d(p_1,f_n(y_i)) +d(f_n(x_j),p_3) +3\eta - c e^{-CR_n}.$$
    Since $\eta\ge d(f_n(y_i),p)$ and $\eta\ge d(p, f_n(x_j))$ we deduce
    $$E(f_n, [-R_n,R_n])\ge d_{13} +\eta - c e^{-CR_n},$$
    which contradicts the energy upper bound if $n$ is large enough, hence proving the claim.

    Thus we have shown that there is either a single transition interval from $q_0 = p_1$ to $q_1 = p_3$ or two transition intervals: one from $q_0 = p_1$ to $q_1 = p_2$ and the other from $q_1 = p_2$ to $q_2 = p_3$.

    We now exclude the first possibility. If there was a single transition, then by choosing $a_n\in I_1$, and going to a subsequence if necessary, we would have $f_n(a_n+\cdot)\to f$ weakly in $H^1_\loc(\R)$ and locally uniformly by compact Sobolev embedding, with $E(f,\R)\le d_{13}$ and $d(f(t),p_1) \le 4\eta$ for  $t$ small enough while $d(f(t),p_3) \le 4\eta$ for  $t$ large enough.  This and the finiteness of $E(f,\R)$ easily implies that $f(-\infty) = p_1$ and $f(+\infty) = p_3$, and thus that $f$ is a heteroclinic joining $p_1$ to $p_3$, contradicting our assumptions on $W$.

    We thus have two transition intervals $I_1 = [y_0,x_1]$ and $I_2 = [y_1,x_2]$ whose complement are $J_0$, $J_1$ and $J_2$. From now on we denote by $I_1 = [\alpha_n,\beta_n]$ and $I_2 = [\alpha_n',\beta_n']$ and we recall that $\alpha_n\ge \frac{-R_n}2$, while $\beta_n'\le \frac{R_n}2$. Moreover  $\alpha_n'-\beta_n$ must converge to $+\infty$ as $n\to +\infty$ as the opposite would imply, as in the previous paragraph, the existence of a heteroclinic joining $p_1$ to $p_3$. Then,  choosing $a_n\in I_1$ and $b_n\in I_2$ and going to a subsequence if necessary, we deduce that $f_n(a_n+\cdot)\to \tilde f$ and $f_n(b_n+\cdot)\to \tilde g$, where $d(\tilde f(t),p_1) \le 4\eta$ for  $t$ small enough and $d(\tilde f(t),p_2) \le 4\eta$ for  $t$ large enough while $d(\tilde g(t),p_2) \le 4\eta$ for  $t$ small enough and $d(\tilde g(t),p_3) \le 4\eta$ for  $t$ large enough. 
    
    Moreover, 
    $$ d_{13} \ge \lim_{n\to +\infty} E(f_n,\R)\ge E(\tilde f, \R) + E(\tilde g,\R).$$
    As above we deduce that 
    $$\tilde f(-\infty) = p_1,\quad \tilde f(+\infty) = p_2,\quad  \tilde g(-\infty) = p_2,\quad\tilde g(+\infty) = p_3,$$
    which implies that  $E(\tilde f,\R) \ge d_{12}$ and $E(\tilde g,\R) = d_{23}$.

    Thus, in view of the upper-bound $d_{13} \ge E(\tilde f, \R) + E(\tilde g,\R),$ we find that $\tilde f$ and $\tilde g$ are heteroclinics connecting, respectively, $p_1$ to $p_2$ and $p_2$ to $p_3$. Shifting $\{a_n\}_n$ and $\{b_n\}_n$ if necessary, we may thus assume that 
    $$ f_n(a_n+\cdot)\to \zeta_{12},\quad f_n(b_n+\cdot)\to \zeta_{23},$$
    in the weak $H^1_\loc$ topology and locally uniformly. Indeed, since $\lim E(f_n)=E(\tilde f)+E(\tilde g)$, it follows that the convergence is, in fact, strong in $H^1_\loc.$
    
    Moreover, from the energy upper-bound we deduce that for any fixed $R>0$ and $n$ large enough depending on $R$, 
    $$ E(f_n,\R\setminus([a_n-R,a_n+R]\cup[b_n-R,b_n+R])) \le \gamma_n + C e^{-cR}.$$

    It remains to prove the quantitative estimates. To this end, we recall \cite{Sch}, Lemma~4.5 which states that there exist $\alpha$ and $\beta$ depending only on $W$ such that if $h:\R\to\R$ satisfies 
    \begin{equation}
    \norm{h-\zeta_{12}}_{H^1(\R)}\le\beta\label{Michelle}
    \end{equation}
    then there exists $t\in\R$ such that 
    \begin{equation}
    \norm{h-\zeta_{12}(\cdot+t)}_{H^1(\R)}^2\leq \alpha\big(E(h)-E(\zeta_{12})\big),\label{Schatzman}
    \end{equation}
    with a similar statement pertaining to $\zeta_{23}.$
    We note that the hypotheses on $W$ in \cite{Sch} differ from ours globally; for example there one assumes the existence of two heteroclinic connections between a pair of wells. However, in a neighborhood of a particular heteroclinic, our situation is identical to that of \cite{Sch} and so the same conclusions hold.
    
    In order to apply this result, we consider extensions of $f_n$ to $\R$ as follows. First, using Lemma~\ref{lem21}, we show that there exists $T_n\in J_1$ such that 
    $$ \|f_n(T_n) - p_2\|^2\le C\big(\gamma_n+e^{-cR_n}\big).$$
  Indeed, Lemma~\ref{lem21} with $I=J_1$ yields 
  $$  E(f_n , J_1)\ge 2\eta +c\left(\min_{t \in[-R_n/2,R_n/2]}|f(t)-p_2|^2+ e^{-CR_n}\right). $$
  Further on $J_i, i=0,2$, we apply Lemma~\ref{lem20} to $f_n$  to obtain, since at one endpoint of $J_0$ (resp $J_2$) $f_n$ is at distance $\eta$ from  $p_1$ (resp. $p_3$),   that $E(f_n , J_i)\ge \eta -c e^{-CR_n}, \ i=0,2$. For the remaining intervals $I_1, I_2$, we note that on $I_1$ , $f_n$ connects values from the geodesic balls centered at $p_1$ to the one centered at $p_2$, while on $I_2$, $f_n$ connects values from the geodesic balls centered at $p_2$ to the one centered at $p_3$, so that 
  $E(f_n , I_1)\ge d_{12} - 2\eta, E(f_n , I_2)\ge d_{23} - 2\eta.$
  Finally using the upper bound $E(f_n , [-R_n,R_n])\le d_{13} +\gamma_n$, we obtain the desired $T_n\in J_1$ with $ \|f_n(T_n) - p_2\|^2\le C\big(\gamma_n+e^{-cR_n}\big)$.

    Next we define $h_n:\R\to\R$ and $h_n':\R\to\R$ as follows. Their respective restriction to $[-R_n/2,T_n]$ and to $[T_n,R_n/2]$ agree with $f_n$. Outside these intervals, we extend $h_n$ to the right of $T_n$ via the energy minimizer connecting $f_n(T_n)$ to $p_2$ and to the left of $-R_n/2$ via the energy minimizer  connecting $f_n(-R_n/2)$ to $p_1$.  Similarly we extend $h_n'$ to the right of  $R_n/2$ via the energy minimizer connecting $f_n(R_n/2)$ to $p_3$ and to the left of $T_n$ via the energy minimizer  connecting $f_n(T_n)$ to $p_2$. 
    
    We now verify that $h_n$ satisfies condition \eqref{Michelle} if $n$ is large enough. First we note that, recalling $I_1 = [\alpha_n,\beta_n]$ and $I_2 = [\alpha_n', \beta_n']$, 
    \begin{equation}\label{sch1}
     \norm{h_n-\zeta_{12}(\cdot-a_n)}_{H^1([\alpha_n,\beta_n])}\to 0\;\text{as}\;n\to\infty.
    \end{equation}
    On the other hand, since 
    $$d(f_n(\alpha_n), p_1)\le \eta,\quad d(f_n(\beta_n), p_2)\le \eta, \quad d(f_n(\alpha_n'), p_2)\le \eta,\quad d(f_n(\beta_n'), p_3)\le \eta,$$
    we have $E(f_n,[\alpha_n,\beta_n])\ge d_{12} - 2\eta$ and $E(f_n,[\alpha_n',\beta_n'])\ge d_{23} - 2\eta$, so that
    \[E(f_n,[-R_n,R_n]\setminus([\alpha_n,\beta_n]\cup[\alpha_n',\beta_n']))\le 4\eta+\gamma_n.\]
    However, on $(-\infty,\alpha_n]$, both $\zeta_{12}(\cdot-a_n)$ and $h_n$ take values in the geodesic ball centered at $p_1$ with radius $4\eta$. Therefore their energy is comparable to the square of their $H^1$ distance to $p_1$ on this interval. Thus, in particular
    \begin{equation}\label{sch2}\|h_n-\zeta_{12}(\cdot-a_n)\|^2_{H^1(-\infty,\alpha_n)}<C\eta,\end{equation}
    for $n$ large enough. Similarly, 
    \begin{equation}\label{sch3}\|h_n-p_2\|^2_{H^1([\beta_n,+\infty])}<C\eta.\end{equation}
    It follows from \eqref{sch1}, \eqref{sch2} and \eqref{sch3} that \eqref{Michelle} holds for $h_n$ if $\eta$ is small enough, depending only on $W$, and for $n$ large enough. The same is true of the $H^1$ distance between $h_n'$ and $\zeta_{23}(\cdot-b_n)$.
 
 In order to exploit the conclusion \eqref{Schatzman}, we will now prove that
 \begin{equation}
 E(h_n)-E(\zeta_{12})<C(\gamma_n +e^{-cR_n})\quad\text{and}\quad
 E(h_n')-E(\zeta_{23})<C(\gamma_n +e^{-cR_n}).\label{stan}
 \end{equation}
 
 We first note that 
 $$E(h_n,\R)=E(f_n,[-R_n/2,T_n])+d(f_n(-R_n/2),p_1)+d(f_n(T_n),p_2)$$
  and 
   $$E(h_n',\R)=E(f_n,[T_n, R_n/2])+d(f_n(R_n/2),p_3)+d(f_n(T_n),p_2).$$
 Next, we have that 
 $d(f_n(-R_n/2),p_1) \le E(f_n,[-R_n, -R_n/2]) + Ce^{-cR_n}$ and $d(f_n(T_n),p_2) \le C(\gamma_n +e^{-cR_n})$. Also, we have that
 $d(f_n(R_n/2),p_3) \le E(f_n,[R_n/2,R_n]) + Ce^{-cR_n}$.
    Thus, we deduce that
    $$E(h_n,\R)\le E(f_n,[-R_n,T_n])+ C(\gamma_n+ e^{-cR_n}),$$
    $$E(h_n',\R)\le E(f_n,[T_n,R_n])+ C(\gamma_n+ e^{-cR_n}).$$
Consequently,
 \[
E(h_n,\R)+E(h_n',\R)\le E(f_n,[-R_n,R_n])+C(\gamma_n +e^{-cR_n}).
 \]
 Moreover, since $E(f_n,[-R_n,R_n])\le d_{12}+d_{23}+\gamma_n$ and $E(h_n,\R)\ge d_{12}$ and 
 $E(h_n',\R)\ge d_{23}$, we deduce that
 \begin{equation}
 E(h_n,\R)-E(\zeta_{12})\le C(\gamma_n +e^{-cR_n})\quad \text{and}\quad
  E(h_n',\R)-E(\zeta_{23})\le C(\gamma_n +e^{-cR_n}).\label{b}
 \end{equation}
 Applying \eqref{Schatzman}, we find that for some $\tilde a_n$ and $\tilde b_n$ one has
$$\|f_n(\tilde a_n+\cdot) - \zeta_{12}(\cdot)\|_{H^1([-R_n/2,T_n])}^2\le C\big(\gamma_n+e^{-cR_n}\big)$$
and
 $$\|f_n(\tilde b_n+\cdot) - \zeta_{23}(\cdot)\|_{H^1([T_n,R_n/2])}^2\le C\big(\gamma_n+e^{-cR_n}\big).$$
This concludes the proof.    
\end{proof}
The following Corollary is a direct consequence of Lemma \ref{Sandier} and the fact that both $|f(T)-p_2|^2$ and $|f(T)-\zeta_{12}(T-a)|^2$ are bounded by $C(\gamma+e^{-cR})$, and hence 
$$|p_2-\zeta_{12}(T-a)|^2\le C(\gamma +e^{-cR}).$$
\begin{cor}\label{cor1} The numbers $a$ and $b$ in the conclusion of Lemma \ref{Sandier}   satisfy
    $$T - a\ge c\min(|\log\gamma|,R),\quad b - T\ge c\min(|\log\gamma|,R),$$
    for some $c>0$ independent of $\gamma$ and $R$.
\end{cor}

\begin{figure}[H]
	\centering
	\includegraphics[width = 0.6\textwidth]{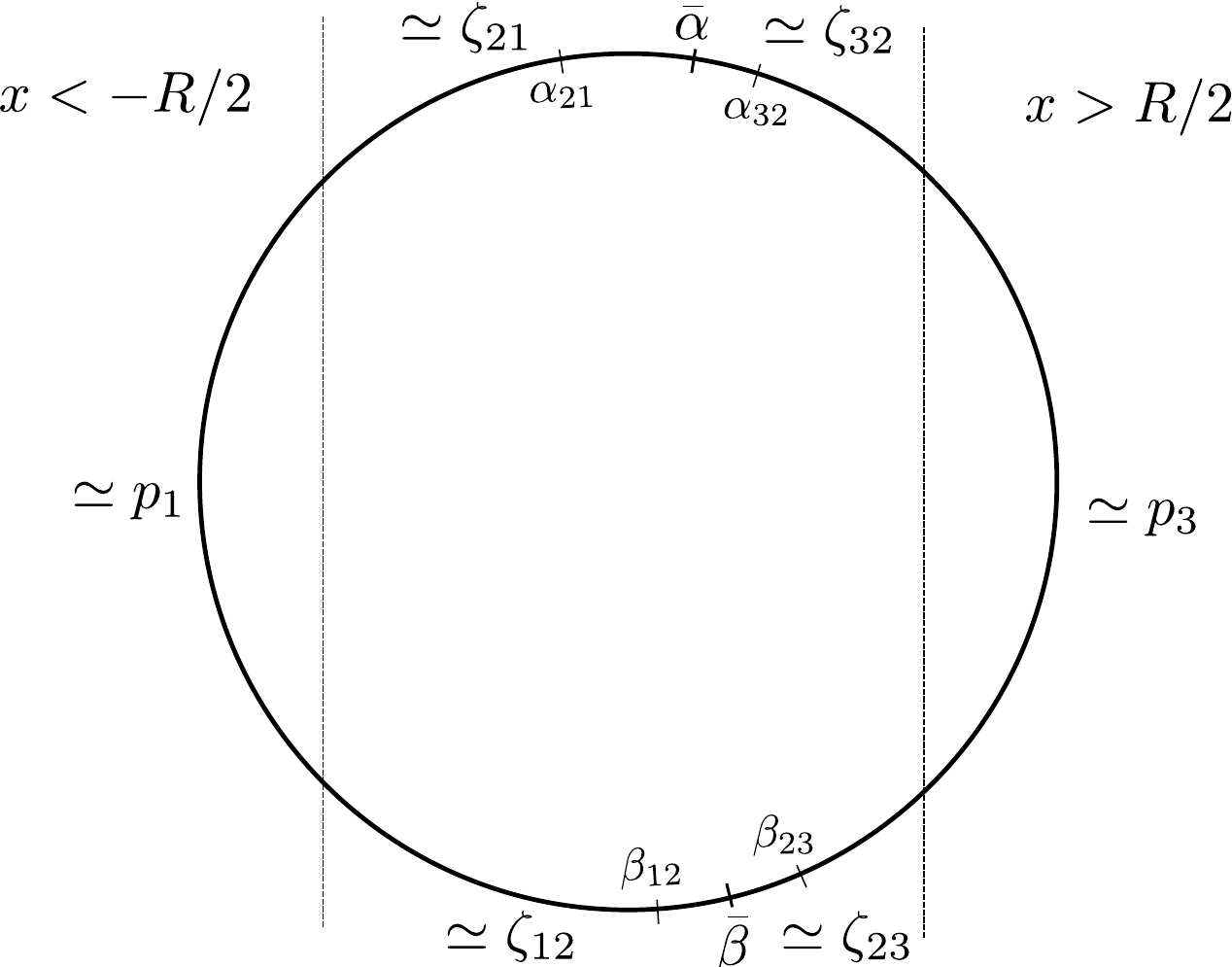}
	\caption{The map $u_{|\partial B_R}$ and the transition angles according to Proposition~\ref{prop1}.}
	\label{circlemap}
\end{figure}

In addition, Lemma \ref{Sandier} also allows us to easily obtain the following Proposition, illustrated in Figure~\ref{circlemap}. Here we denote by $s$ the arclength variable on the circle $\pt B_R$ such that $s(Re^{i\theta}) = R\theta$ for $\theta\in (-\pi,\pi]$.

\begin{pro}\label{prop1} There exists constants $C$, $\eta>0$ and $R_0$, $\gamma_0>0$ such that for any $R>R_0$ and any $0<\gamma<\gamma_0$, if $u:\pt B_R\to \R^2$ is such that 
\begin{gather}
\text{$|u(x,y) - p_1| < \eta$ for any $(x,y)\in\pt B_R$ with $x< -R/2$, }\label{p1}\\
\text{$|u(x,y) - p_3| < \eta$ for any $(x,y)\in\pt B_R$ with $x>R/2$, }\label{p3}\\
E(u,\pt B_R)\le 2 d_{13} +\gamma,\label{up}
\end{gather}
then there exists angles $\pi/4<\alpha_{32}< \ba <\alpha_{21}<3\pi/4$ and $-3\pi/4<\beta_{12}<\bb<\beta_{23}<-\pi/4$ such that the following holds. Writing $I_{21} = [\ba,3\pi/4]$, $I_{32} = [\pi/4,\ba]$ and $I_{12} = [-3\pi/4,\bb]$, $I_{23} = [\bb,-\pi/4]$ , we have 
$$\left\|u(Re^{i(\alpha_{21} +\cdot/R)}) - \zeta_{21}(\cdot)\right\|_{H^1(I_{21})}^2 < C\gamma,\quad
\left\|u(Re^{i(\alpha_{32} +\cdot/R)}) - \zeta_{32}(\cdot)\right\|_{H^1(I_{32})}^2 < C\gamma,\quad $$
\begin{equation}
\left\|u(Re^{i(\beta_{12} +\cdot/R)}) - \zeta_{12}(\cdot)\right\|_{H^1(I_{12})}^2 < C\gamma,\quad
\left\|u(Re^{i(\beta_{23} +\cdot/R)}) - \zeta_{23}(\cdot)\right\|_{H^1(I_{23})}^2 < C\gamma.
\label{yvonne}\end{equation}
and such that for any $s\in [-R\pi,R\pi]\setminus\cup_{jk} I_{jk}$,
\begin{equation}
d(u(Re^{is/R}),P) < \eta.
\label{sylvie}
\end{equation}
Finally, 
\begin{equation}
R\min(\ba-\alpha_{32},\alpha_{21}-\ba, \beta_{23}-\bb,\bb-\beta_{12})\ge c\min(|\log\gamma|,R).\label{minab}
\end{equation}
\end{pro}

\section{Upper bound}

Here we establish a key upper bound:

\begin{pro}\label{prop4} There exists $C>0$ such that for any  $R>C$ and $1/C>\eta>0$, $1/C>\gamma>0$ the following holds.
    
    If $u:\partial B_R\to\R^2$ satisfies the hypotheses of Proposition~\ref{prop1} and is such that $u(x,y)$ is $\eta$-close to $p_1$ if $x<-\eta R$ and $\eta$-close to $p_3$ if $x>\eta R$, then  for any $C\le \sigma$, $\rho \le R/C$ it holds that
\beq E(u,B_R)\le d_{12} L_1  + d_{23} L_3 + C\(\frac{\rho^2}\sigma+\frac\sigma\rho e^{-\ell}+Re^{-(\ell+\rho)}+ \eta\ell e^{-c\ell}+e^{-c\ell} +\gamma +\eta^2\),\label{upp1}
\eeq
where, using the angles $\alpha_{32}$, $\alpha_{21}$, $\beta_{12}$, $\beta_{23}$, $\ba$, $\bb$ from the conclusion of Proposition~\ref{prop1}, we have used the notation
\beq L_1 = R|e^{i\alpha_{21}} - e^{i\beta_{12}}|,\; L_3 = R|e^{i\alpha_{32}} - e^{i\beta_{23}}|,\; \ell = \frac{1}{2}R\,\alpha_{min}
\label{lengths}
\eeq
and
\[
\alpha_{min}:=\min(\ba-\alpha_{32},\alpha_{21}-\ba, \beta_{23}-\bb,\bb-\beta_{12}).
\]
\end{pro}

\begin{proof}[Proof of Proposition~\ref{prop4}]
    The Proposition is proved by constructing a comparison map $v: B_R\to\R^2$ coinciding with $u$ on $\partial B_R$. 

\noindent [{\bf Step 1}. Annulus] We begin by defining a competitor in the annulus $A = B_R\setminus B_{R-1}$. To this end, it will be convenient to define the following modification of the heteroclinics $\zeta_{ij}$ using the parameter $\ell$ from \eqref{lengths}. We set
\begin{equation}\zeta_{ij}^{\ell}(s) =
     \begin{cases} 
        p_i & \text{if $s\le -\ell$,}\\
        \zeta_{ij}(s) & \text{if $|s|<\ell - 1$,}\\
        p_j & \text{if $s\ge \ell$,}
    \end{cases}\label{zetaell}
\end{equation}
    and we take $\zeta_{ij}^\ell$ to be a linear interpolation on the intermediate intervals.  It is easy to check that 
    \begin{equation}\label{newlapprox} E(\zeta_{ij}^{\ell}, [-\ell,\ell])\le d(p_i,p_j) + e^{-c\ell},\quad \|\zeta_{ij}^{\ell}-\zeta_{ij}\|_{L^2(\R)}^2\le e^{-c\ell}.\end{equation}

Adopting the notation from Proposition~\ref{prop1}, we note that the line through $e^{i\bar{\alpha}}$ and $e^{i\bar{\beta}}$, which we denote by $L$, separates $B_R$ into two regions, one, say $B_R^{{\rm left}}$, containing the line segment joining $Re^{i\alpha_{21}}$ and $Re^{i\beta_{12}}$ and $B_R^{{\rm right}}$ containing the line segment joining $Re^{i\alpha_{32}}$ and $Re^{i\beta_{23}}$.

Denoting by $S_{12}$ the first of these two segments and by $S_{23}$ the second, we define a competitor, say $v$, in $A$ as the linear interpolation in the radial direction between the values of $u$ on $\partial B_R$  and the function, say $V:\partial B_{R-1}\to\R^2$ given by
\begin{equation}
V(x,y):=\left\{\begin{matrix} \zeta^\ell_{12}({\rm dist}\big((x,y),S_{12}\big)&\text{for}\;(x,y)\in\partial B_{R-1}\cap B_R^{{\rm left}}\\
&\\
 \zeta^\ell_{23}({\rm dist}\big((x,y),S_{23}\big)&\text{for}\;(x,y)\in\partial B_{R-1}\cap B_R^{{\rm right}}.\end{matrix}\right.\label{Vdef}
\end{equation}
Here ${\rm dist}\big((x,y),S_{12}\big)$ denotes the signed distance function to $S_{12}$, taken to be negative to the left of $S_{12}$ and 
${\rm dist}\big((x,y),\ell_{23}\big)$ is the signed distance function to $S_{23}$, taken to  be negative to the left of $S_{23}$. We note that, so defined, the function $V$ is continuous taking value $p_2$ at the two points $\partial B_{R-1}\cap\partial B_R^{{\rm left}}\cap \partial B_R^{{\rm right}}.$ We also note that in light of the definition of $\ell$, cf. \eqref{lengths},  the segments $S_{12}$ and $S_{23}$ are both more than $\ell$ distance from the line $L$.

 Estimating the energetic cost of $v$ in the annulus $A$ is entirely analogous to a similar calculation carried out \cite{EP}, step 3 of the proof of Theorem 1.3. Therefore, we omit the details. One finds that with two transitions from $p_1$ to $p_2$ and two transitions from $p_2$ to $p_3$, one has
$$
E(v,A)\le 2\(d_{12} + d_{23}\) + C\(e^{-c \ell}+\gamma+\eta^2\),
$$
through an appeal to \eqref{newlapprox} and \eqref{yvonne}. Here we have used the fact that the direction of the segments $S_{12}$ and $S_{23}$ and the radial direction at their endpoints make angles bounded by $C\eta$ resulting in the $O(\eta^2)$ term above. 

\noindent[{\bf Step 2}. Defining $v$ in `triangles' inside $B_{R-1}$ bordering $\partial B_{R-1}$]   The line $L$ passing through $e^{i\ba}$ and $e^{i\bb}$ separates $B_{R-1}$ into two parts: one to the left, which we denote by $D$,  containing the segment $S_{12}\cap B_{R-1}$, and the other one containing the segment $S_{23}\cap B_{R-1}$. We then introduce a coordinate system $(x',y')$ to define the competitor $v$ in the region $D$ by taking the segment $S_{12}$ to coincide with the $y'$ axis, with 
$S_{23}\cap B_{R-1}$ extending in $y'$ coordinates from say $y'=y_-$ to $y'=y_+$.
 We also take 
 the $x'$-axis to go through the middle of $S_{23}\cap B_{R-1}$ so that $y_-=-y_+$. Before proceeding, we remark that in these coordinates, the function $V$
defined in \eqref{Vdef} is simply a function of $x'$ for $(x,y)\in\partial B_{R-1}\cap B_R^{{\rm left}}$.

 In these coordinates, consider the strip $\{(x',y')\mid |x'| < \ell\}$. Its intersection with $B_{R-1}$ may be decomposed as a rectangle $\rect$, defined as the largest rectangle of the form $\{(x',y'):\, |x'| < \ell,\ y'\in[y_-,y+]\}$ inscribed in $B_{R-1}$, and  two ``triangles'' $C_+$ and $C_-$, each having one side which is in fact an arc in $\partial B_{R-1}$. The other two sides are parallel to the $x'$ and $y'$ axis, respectively. Since $Re^{i\alpha_{21}}$ and $Re^{i\beta_{12}}$ both must lie to the right of the line $x=-\eta R$ by assumption, it follows that the ends of the segment $S_{12}\cap B_{R-1}$ are at a distance of at most $\eta R$ from the $y$-axis, and so the ratio of the length of the $y'$-side to the $x'$-side for both of the aforementioned `triangles' is bounded by $C\eta$. Consequently, the length of the $y'$ side of these triangles is $O(\eta\ell)$.

    On $C_+$ and $C_-$, we define $v$ to be independent of $y$  so that 
    $$v(x',y_+)  = v(x',-y_+) = \zeta_{12}^{\ell}(x'),$$
    and 
    \begin{equation}
    \begin{split}
    E(v,C_+\cup C_-)  
    &\le  E(\zeta_{12}^{\ell})\times\( \left|(R-1)e^{i\alpha_{21}} -(R-1) e^{i\beta_{12}}\right| - 2y_+\)\\
    &\le d_{12} \( \left|(R-1)e^{i\alpha_{21}} - (R-1)e^{i\beta_{12}}\right| - 2y_+\)+ C\eta\ell e^{-c\ell}.
    \end{split}\label{split}
\end{equation}

\noindent[{\bf Step  3}, Dilation] Though tempting, it is too costly energetically to simply take our competitor $v$ to be $\zeta_{12}^{\ell}(x')$ throughout $D$. Therefore, we must further separate the center of this modified heteroclinic from the line $L$. To this end,
we proceed to define $v$ first on 
    $$D_+ := D\cap \R\times[y_+-\sigma,y_+],\quad D_- := D\cap \R\times[-y_+, \sigma-y_+].$$ 
    By enforcing the condition $v(x',y') = v(x',-y')$ on $D_+\cup D_-$, in fact, we only need to work on $D_+$.
    
    The idea is to interpolate between $\zeta_{12}^{\ell}(x')$ and $\zeta_{12}^{\ell+\rho}(x'+\rho)$, in the $y'$ variable for a parameter $\rho$ to be chosen judiciously later. That is, for $(x',y')\in D_+$ we take our competitor $v$ to be given by
    \begin{equation}
    v(x',y')= \zeta_{12}^{\ell+\frac{\rho}{\sigma}(y_+-y')}\big(x'+\frac{\rho}{\sigma}(y_+-y')\big).\label{vinterp}
    \end{equation}
    In light of \eqref{newlapprox}, the energy due to the $x'$ derivative and the potential term in $D_+$ is bounded by
    \begin{equation}
    \int_0^\sigma d_{12} + e^{-c(\ell+t\rho/\sigma)}\,dt
     \le  \sigma d_{12}+C \frac{\sigma}\rho e^{-c\ell}.\label{xprime}
     \end{equation}
      From \eqref{vinterp} and the definition of $\zeta_{12}^\ell$, it is straightforward to deduce that 
     \begin{equation}\frac12\int_{D_+}\left|\frac{\partial v}{\partial y'} \right|^2 \le \( \frac\rho\sigma \)^2 \int_0^\sigma  d_{12} + e^{-c(\ell+t\rho/\sigma)}\,dt\leq C\frac{\rho^2}{\sigma}.\label{yprime}
     \end{equation}
      
    As noted at the beginning of this step, the same bounds hold for $D_-$.

\noindent[{\bf Step  4},  Remainder and conclusion] It remains to define $v$ in $D\cap \R\times[y_-+\sigma,y_+-\sigma].$ We simply require $v$ to be independent of $y'$ there, so that $v=v(x')=\zeta_{12}^{\ell+\rho}(x'+\rho).$ Hence, in this region the energy bound takes the form 
    
    \begin{equation}E(v,D\cap \R\times[y_-+\sigma,y_+-\sigma]) \le (y_+-y_- - 2\sigma) (d_{12} + e^{-c(\ell+\rho)}).\label{bulkD}
    \end{equation}
    
    Adding the upper bounds \eqref{split}, \eqref{xprime}, \eqref{yprime} and \eqref{bulkD},  recalling that $y_- = -y_+$, and noting that $y_+<R$, we obtain
   \begin{multline*}E(v,D)\le (2y_+- 2\sigma) \(d_{12} + e^{-c(\ell+\rho)}\)+ C\frac{\rho^2}\sigma +2\sigma d_{12}\\+C \frac{\sigma}\rho e^{-c\ell} + d_{12} \( \left|(R-1)e^{i\alpha_{21}} - (R-1)e^{i\beta_{12}}\right| - 2y_+\)+ C\eta\ell e^{-c\ell}.
   \end{multline*}
    Therefore,
    $$E(v,D)\le d_{12} \left|(R-1)e^{i\alpha_{21}} - (R-1)e^{i\beta_{12}}\right| +
    C\(R  e^{-c(\ell+\rho)} + \frac{\rho^2}\sigma +  \frac{\sigma}\rho e^{-c\ell} + \eta\ell e^{-c\ell}\).$$
    Carrying out an analogous construction in the portion of $B_{R-1}$ lying to the right of the line $L$, that is, with $\zeta_{12}^\ell$ replaced by $\zeta_{23}^\ell$, etc, we deduce that 
    \begin{multline*}E(v,B_{R-1})\le d_{12} \left|(R-1)e^{i\alpha_{21}} - (R-1)e^{i\beta_{12}}\right| +
    d_{23} \left|(R-1)e^{i\alpha_{32}} - (R-1)e^{i\beta_{23}}\right|+\\ + C\(R  e^{-c(\ell+\rho)} + \frac{\rho^2}\sigma +  \frac{\sigma}\rho e^{-c\ell} + \eta\ell e^{-c\ell}\).
    \end{multline*}
    Then because all of the transition angles are within $C\eta$ of being $\pm\pi/2$ due to the hypothesis, we have
    $$ \left|Re^{i\alpha_{21}} - Re^{i\beta_{12}}\right| = \left|(R-1)e^{i\alpha_{21}} - (R-1)e^{i\beta_{12}}\right| + 2+ O(\eta^2),$$
    with a similar statement holding if we replace $\alpha_{21}$ and $\beta_{12}$ by $\alpha_{32}$ and $\beta_{23} $. Thus, once we add in the cost of the annulus interpolation  $E(v,A)$ from Step 1, we finally arrive at the upper bound 
    $$E(v,B_R)\le d_{12} L_{12} + d_{23} L_{23} + C\(R  e^{-c(\ell+\rho)} + \frac{\rho^2}\sigma +  \frac{\sigma}\rho e^{-c\ell} + \eta\ell e^{-c\ell}+e^{-c\ell} +\gamma +\eta^2\).$$

\end{proof}

\begin{lemma}\label{cvh13}
    Assume $u:\R^2\to\R^2$ is a minimizer of $E$ and that there exists a sequence $R_n\to +\infty$ such that, possibly after rotating the coordinates, $u(R_n\;\cdot)\to H_{13}$ in $L^1$. Then {\em any} sequence $R_n\to +\infty$ is such that, possibly after rotating the coordinates, $u(R_n\;\cdot)\to H_{13}$.
\end{lemma}

\begin{proof} Assume by contradiction that no subsequence of $u(R_n\;\cdot)$ converges to $H_{13}$ modulo a rotation, then a subsequence must converge to $H_{12}$ or $H_{23}$ modulo a rotation, which implies,  from \cite{EP}, Theorem~1.3, that, modulo a rotation, $u(x,y) = \zeta_{12}(x)$ or $u(x,y) = \zeta_{23}(x)$, which contradicts the assumptions.
\end{proof}

\begin{lemma} \label{lowexp}
    Under the same hypothesis, there exist $R_0>0$ such that for any $R>R_0$, 
    $$E(u,\pt B_R)\ge 2 d_{13} - e^{-CR}.$$
\end{lemma}

\begin{proof} It suffices to prove that for any sequence $R_n\to +\infty$ the bound is satisfied with a constant $C$ independent of $n$. We choose $\eta>0$ small enough depending on $W$, to be determined below.
    
    Since $u(R_n\;\cdot)$ converges in $L^1$ to $H_{13}$ modulo a rotation, for any $R$ large enough in the sequence there exists 
    $\alpha$, $\alpha'$ and $\beta$, $\beta'$ with
    $$\pi/2-\eta<\alpha<\pi/2<\alpha'<\pi/2+\eta, 
    \quad -\pi/2-\eta<\beta<-\pi/2<\beta'<-\pi/2+\eta$$
   such that $u(Re^{i\alpha})$, $u(Re^{i\beta})$ are $\eta$-close to $p_1$ and such that $u(Re^{i\alpha'})$, $u(Re^{i\beta'})$ are $\eta$-close to $p_3$.

    
    To bound from below the energy of $u$ on $\partial B_R^+ = \{Re^{i\theta}\mid\theta\in[0,\pi]\}$ we let $v(s) = u(Re^{is/R}),$ so that $E(u,\partial B_R^+) = E(v,[0,\pi R]).$ Then 
    $$E(v,[0,\pi R])  = E(v,[0,R\alpha'])+E(v,[R\alpha', R\alpha])+E(v,[R\alpha,R\pi]).$$
    We have $E(v,[R\alpha', R\alpha])\ge d(v(R\alpha'), v(R\alpha))$ and, using Lemma~\ref{lem20}, we find  that 
    $E(v,[0,R\alpha'])\ge d(v(R\alpha'), p_3)- Ce^{-cR}$ while $E(v,[R\alpha,R\pi]) \ge d(v(R\alpha), p_1)- Ce^{-cR}$. Therefore
    $$E(u,\partial B_R^+) \ge  d_{13} - e^{-cR}.$$
    The same holds on $\partial B_R^-$, proving the result.
\end{proof}

We are now ready to prove the following important upper bound:

\begin{pro}\label{prop2} Assume $u:\R^2\to\R^2$ is a minimizer of $E$ such that $u(R\;\cdot)\to H_{13}$, possibly after taking subsequences and modulo a rotation of the coordinates. 
    
    Then there exists  sequences of positive numbers  $R_n\to +\infty$ and $\eta_n\to 0$, $\gamma_n\to 0$, such that the restriction of $u$ to $\pt B_{R_n}$ satisfies the hypothesis of Proposition~\ref{prop1} with the parameters $\eta_n$, $\gamma_n$. Moreover, denoting $\alpha_{21}$, $\alpha_{32}$, $\beta_{12}$, $\beta_{23}$ the transition angles (which depend on $n$) given by Proposition~\ref{prop1} we have, as $n\to +\infty$, 
    $$E(u,B_{R_n})\le R_n\(|e^{i\alpha_{21}}- e^{i\beta_{12}}| d_{12} + |e^{i\alpha{32}}- e^{i\beta_{23}}| d_{23}\) + o(1).$$
\end{pro}


\begin{proof}[Proof of Proposition~\ref{prop2}]
    From Lemma~\ref{cvh13}, any sequence $R_n\to +\infty$ has a subsequence such that, modulo a rotation, $u(R_n\;\cdot)$ converges to $H_{13}$. It follows that, given $\eta$, $\gamma>0$, if $n$ is large enough then \eqref{p1}, \eqref{p3} are satisfied for every $R\in[\eta R_n,R_n]$. Moreover, from Lemma~\ref{lowexp}, for $R\in[\eta R_n,R_n]$ we have 
    $$E(u,\pt B_R) \ge 2 d_{13} - e^{-cR}.$$

    On the other hand Proposition~\ref{gamma} implies, after rescaling, that  as $n\to +\infty$,
    $$E(u,B_{R_n}) \le 2R_nd_{13}+ o(R_n),$$
    since the diameter of $B_{R_n}$ is $2R_n$.

    Combining the two, we find that \eqref{up} is satisfied with $\gamma_n = o(1)$, for most $R\in[\eta R_n,R_n]$. Hence \eqref{p1}, \eqref{p3} and \eqref{up} are satisfied for some sequence $\{R_n\}$ converging to $+\infty$ with $\gamma_n = o(1)$.

    We then apply Proposition~\ref{prop4} with $\sigma = \rho = R^\theta$ for some fixed $0<\theta<1$, for any $R$ in our sequence. We find that all the error terms may be absorbed in $CR^\theta$ so that  
    $$E(u,B_R) \le d_{12} L_1  + d_{23} L_3 + CR^\theta,$$
    Moreover, from Lemma~\ref{lowexp}, we have for any $r>R_0$ that 
    $$E(u,\pt B_r)\ge 2 d_{13} - e^{-Cr}.$$
    It follows that 
    \begin{equation}
        \begin{split}
            \int_{R/2}^R E(u,\pt B_r) -  2 d_{13} \,dr &= E(u,B_R) - E(u,B_{R/2}) - R d_{13}\\
            &\le E(u,B_R) - \int_{R_0}^{R/2} 2 d_{13} - e^{-Cr}\,dr - R d_{13}\\
            &\le E(u,B_R) - 2R d_{13} + C \\
            &\le L_1 d_{12} + L_3 d_{23} + CR^\theta - 2Rd_{13} \\
           &\le  { CR^\theta,}
        \end{split}
    \end{equation}
    {since} $L_1$, $L_3\le 2R$.
    { Hence} there must exist $R'\in[R/2,R]$ such that 
    $$E(u,\pt B_{R'})\le 2d_{13} + CR^{\theta-1}.$$

    We denote $\{{R'}_n\}$ a subsequence of $\{{R}'\}$ such that $u({R'}_n\cdot)\to H_{13}$ modulo a rotation. We may then apply Proposition~\ref{prop1} on $\partial B_{R'}$, with $\eta_n$ tending to $0$ and $\gamma = C{({R'}_n})^{\theta -1}$, which also tends to $0$ as $n\to +\infty$. It yields transition angles $\alpha_{ij}$, $\beta_{ij}$ such that 
    $$c\log(({{R'}_n})^\theta)\le \ell := \frac{{R'}_n}2\min(|e^{i\alpha_{21}} - e^{i\alpha_{32}}|,|e^{i\beta_{12}} - e^{i\beta_{23}}|).$$
    Therefore, Proposition~\ref{prop4} applied with $\rho = {({R'}_n})^\e$ and $\sigma = {({R'}_n})^{3\e}$, with $\e>0$ chosen small enough (taking $\e < c\theta$ works), yields
    $$E(u,B_{{R'}_n})\le L_1 d_{12} + L_3 d_{23} + o(1).$$
\end{proof}

\section{Proof of Theorem~\ref{main}}

We now proceed to prove our main result, by contradiction. It is proved in \cite{EP} that if $u$ is a minimizer of $E$ and $\{R_n\}$ is a sequence of radii tending to $+\infty$ such that $u(R_n\cdot)$ converges to $H_{12}$ in $L^1$, then $u(x,y) = \zeta_{12}(x)$. The same holds if $\{1,2\}$ is replaced by $\{2,3\}$. In light of Proposition~\ref{gamma}, it thus remains to obtain a contradiction under the assumption that, for any $\{R_n\}$ converging to $+\infty$, there exists a subsequence such that, possibly after a rotation of the coordinates, $\|u(R_n\;\cdot) - H_{13}(R_n\;\cdot)\|_{L^1}$ converges to $0$ as $n\to +\infty$, which we assume henceforth. 
 
The contradiction will be obtained by defining $\{R_n\}$, $\{\eta_n\}$, $\{\gamma_n\}$ to be as in the conclusion of Proposition~\ref{prop2}. Our aim is to derive a lower bound for $E(u,B_{R_n})$ that will contradict the upper bound of Proposition~\ref{prop2}, for large enough $n$.

The map  $u(R_n\cdot)$ converges to $H_{13}$ in $L^1$, and therefore locally uniformly outside the $y$-axis, see remark ~\ref{locunif}. Therefore, for any $\eta>0$ we may choose $n_0$ large enough so that $d(u(x,y), p_1)<\eta$ for any $(x,y) \in B_{R_{n_0}}$ such that $x<-\eta R_{n_0},$ while $d(u(x,y), p_3)<\eta$  if $(x,y)\in B_{R_{n_0}} $ is such that  $x>\eta R_{n_0}$. Relabeling the sequence we let $n_0 = 0$ and $R_{n_0} = R_0$. We thus have
\begin{equation}
    \parbox{.9\textwidth}{$d(u(x,y), p_1)<\eta$ in $\mathcal R_\lft = B_{R_0}\cap\{x<-\eta R_0\}$ and $d(u(x,y), p_3)<\eta$ in $\mathcal R_\rght = B_{R_0}\cap\{x>\eta R_0\}$.}\label{central}
\end{equation}

We now proceed to prove the desired lower-bound for $E(u,B_{R_n})$. First we extend $u$ to $\tilde u$ defined on $B_{R_n+1}$, which will allow us to clean up the boundary data. We claim that we may define $\tilde u$ on the annulus $A_n = B_{{R_n}+1}\setminus B_{R_n}$ so that, as $n\to +\infty$
\begin{equation}
    E(\tilde u, A_n) \le 2(d_{12}+d_{23}) + o(1),\label{annular}
\end{equation}
and so that $\tilde u$ restricted to $\partial B_{R_n+1}$ is equal to $\zeta^\ell_{ij}$, defined by \eqref{zetaell},  on an arc of length $\ell$ centered at the angle $\alpha_{ij}$ or $\beta_{ij}$ and $\tilde u$ is constant equal to one of the wells on each of the connected components of the complement of these arcs.   
Here the angles $\alpha_{ij}$, $\beta_{ij}$ are the transition angles (which depend on $n$) given by Proposition~\ref{prop1}. Note that $\tilde u$ depends on $n$ as well, but only on the annulus $A_n$. The construction of $\tilde u$ is by radial interpolation, as in the proof of Proposition~\ref{prop4}, Step~1, and the estimate \eqref{annular} follows from the same arguments, hence is omitted.

Again using the fact that  $u$ converges to $H_{13}$ locally uniformly to $p_1$ on the half-plane $x<0$ and to $p_3$ on the half-plane $x>0$, we have that  transition angles tend to $\pm\pi/2$ as $n\to +\infty$. Therefore
$$ \left|(R_n+1)e^{i\alpha_{21}} - (R_n+1)e^{i\beta_{12}}\right| = \left|R_ne^{i\alpha_{21}} - R_ne^{i\beta_{12}}\right|+ 2 + o(1), $$
and a similar statement holds if we replace $\alpha_{21}$ and $\beta_{12}$ by $\alpha_{32}$ and $\beta_{23} $. Thus, from \eqref{annular} and the upper-bound of Proposition~\ref{prop2} we deduce that 
$$E(\tilde u, B_{R_n+1}) \le (R_n+1)\(|e^{i\alpha_{21}}- e^{i\beta_{12}}| d_{12} + |e^{i\alpha_{32}}- e^{i\beta_{23}}| d_{23}\) + o(1).$$
From now on we drop the tilde and denote $u$ the extension of $u$ to $B_{R_n+1}$, and write $R$ instead of $R_n+1$ for simplicity. The above upper-bound thus becomes 
\begin{equation}
    E(u, B_{R}) \le R\(|e^{i\alpha_{21}}- e^{i\beta_{12}}| d_{12} + |e^{i\alpha_{32}}- e^{i\beta_{23}}| d_{23}\) + o(1).\label{up2}
\end{equation}
We proceed to prove a lower-bound contradicting \eqref{up2} for large enough $n$.

Let 
$$a_1 = Re^{i\alpha_{21}},\quad b_1 = Re^{i\beta_{12}},\quad a_3 = Re^{i\alpha_{32}},\quad b_3 = Re^{i\beta_{23}}.$$
We define $\zeta:B_R\to\R$ such that 
\begin{equation}\label{Fred}\|\nabla\zeta\|_\infty\le 1,\quad \zeta(a_1) - \zeta(b_1) = |a_1-b_1|,\quad \zeta(a_3) - \zeta(b_3) = |a_3 - b_3|.
\end{equation}
Such a $\zeta$ exists: indeed, since $$\pi>\alpha_{21} > \alpha_{32} > \beta_{23}>\beta_{12} > -\pi,$$
the minimal connection joining the points $a_1$, $a_3$ to the points $b_1$, $b_3$ is exactly the union of the segments $[a_1,b_1]$ and $[a_3,b_3]$, see \cite{BCL}. An explicit formula for $\zeta$ is
\begin{equation}\label{zeta}\zeta = \max\(\zeta_1,\zeta_3\),\quad \zeta_1(x) = x\cdot  u_1 - \lambda_1,\quad \zeta_3(x) = x\cdot u_3 - \lambda_3,\end{equation}
where $u_1 = (a_1-b_1)/|a_1-b_1|$ and $ u_3 = (a_3-b_3)/|a_3-b_3|$, and $\lambda_1$, $\lambda_3$ are suitably chosen real numbers such that $\zeta$ satisfies the required properties. 
\begin{figure}[H]
	\centering
	\includegraphics[width = 0.6\textwidth]{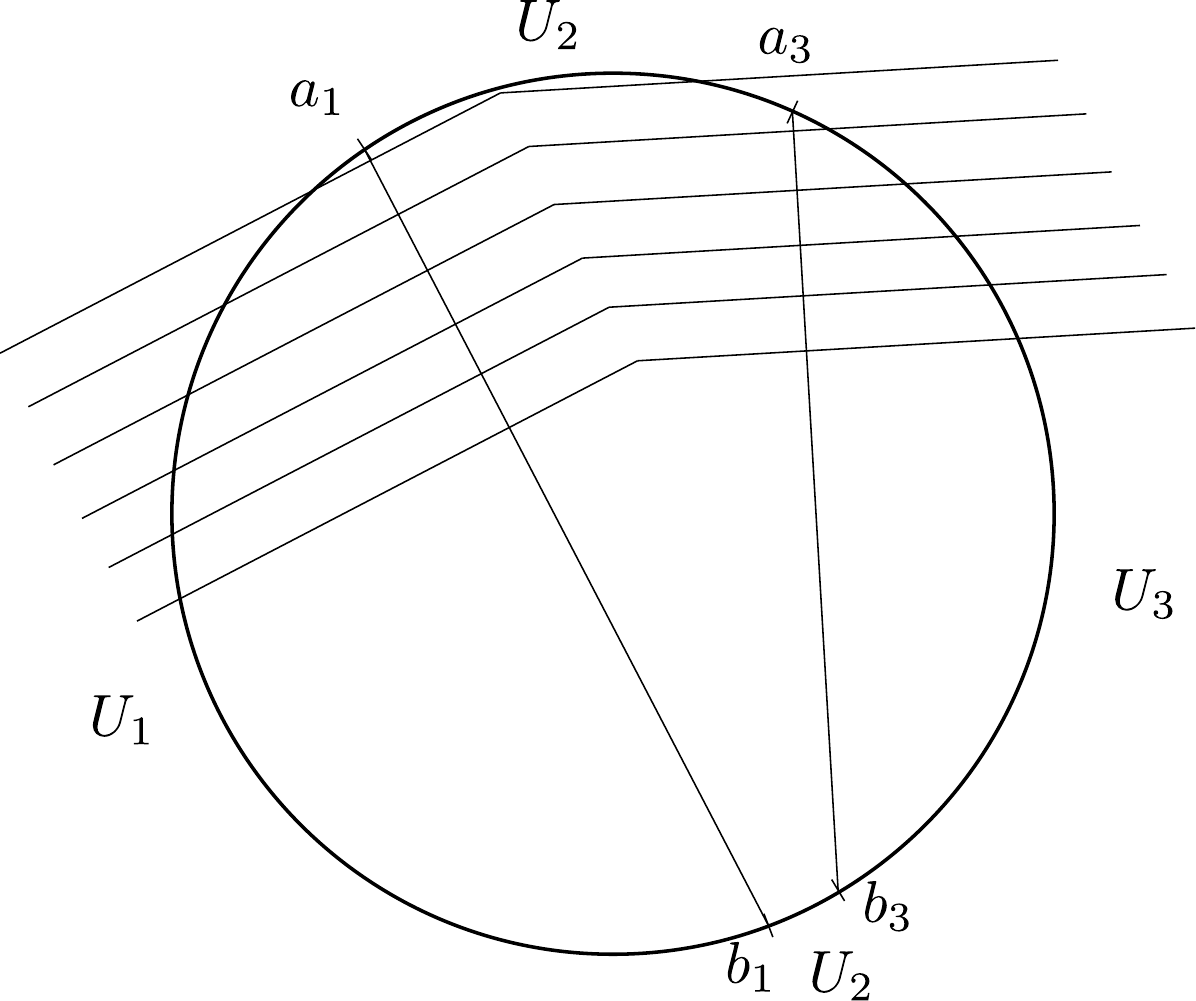}
	\caption{The disk $B_R$ and some level-sets of the function $\zeta$.}
	\label{slicing}
\end{figure}

Since $\zeta_1$ and $\zeta_3$ are linear functions each of whose level sets (see figure~\ref{slicing})
necessarily consist of  a line segment with endpoints on $\partial B_R$, it follows that each nonempty level set of $\zeta$ consists of one or two segments which can be consecutive or disjoint, hence intersect $\partial B_R$ at either $2$ or $4$ points.  The desired lower-bound will be obtained by integrating w.r.t the variable $t$ the lower-bound we next compute on the level set $\ga_t = \{\zeta = t\}$.

The essential parameter in the lower-bound is the boundary condition at the points of  $\ga_t\cap\partial B_R$. The domain $\overline B_R\setminus([a_1,b_1]\cup [a_3,b_3])$ has three connected components. We denote by  $U_1$, $U_2$, $U_3$ the intersections of these three components with $\partial B_R$, numbered so that $u$ restricted to $U_i$ is closest to the well $p_i$, in a way that we will quantify below. We now examine to which $U_i$'s the points in $\ga_t\cap\partial B_R$ belong to.
    
    To this end we note that on $[a_1,b_1]$, the function $\zeta$ is one-to-one, since $\zeta(a_1) - \zeta(b_1) = |a_1 - b_1|$ and $\|\nabla\zeta\|_\infty\le 1$. For the same reason $\zeta$ restricted to $[a_3,b_3]$ is one-to-one, and therefore $\ga_t$ intersects either of these segments  at most once. Four cases may thus occur, depending on the value of $t$:
    \begin{itemize}
        \item[a)] $t\in[\zeta(a_1),\zeta(b_1)]$ and $t\notin[\zeta(a_3),\zeta(b_3)]$. Then $\ga_t$ intersects $[a_1,b_1]$ only. Thus $\ga_t$ contains a segment or broken line which connects $U_1$ to $U_2$.
        \item[b)] $t\notin[\zeta(a_1),\zeta(b_1)]$ and $t\in[\zeta(a_3),\zeta(b_3)]$. Then $\ga_t$ contains a segment or broken line which connects $U_2$ to $U_3$.
        \item[c)] $t\in[\zeta(a_i),\zeta(b_i)]$, for both $i=1,3$. Then $\ga_t$ contains either a segment or broken line connecting  $U_1$ to $U_3$, or two segments: one connecting $U_1$ to $U_2$ and the other connecting $U_2$ to $U_3$
        \item[d)] $t\notin[\zeta(a_1),\zeta(b_1)]$ and $t\notin[\zeta(a_3),\zeta(b_3)]$. Then every connected component of $\ga_t$ starts and ends in the same $U_i$. We will use the trivial lower-bound $E(u,\ga_t)\ge 0$ in this case.
    \end{itemize}
    
    We now quantify the distance of $u(q)$, where  $q\in \ga_t\cap U_i$, to the well $p_i$. Let $t=\zeta_1(a_1) - \delta$. Then $q = a_1 e^{i\alpha_1}$ for some $\alpha_1$ and, if we denote $z = a_1 - \delta u_1$, then $q-z$ is orthogonal to $ u_1$. Therefore, writing $a_1 = R u_1 e^{i\theta_1}$ and recalling that the inner product is given by $\Re (\bar v w)$, we  have $R e^{i(\alpha_1+\theta_1)} -  R e^{-i\theta_1} + \delta\in i\R$, i.e. 
    \begin{equation}R \cos(\alpha_1+\theta_1) -  R \cos \theta_1 + \delta=0.\label{ortho}\end{equation}

    From \eqref{ortho} we deduce that either $|\alpha_1|>R^{-1/2}$ or, using the mean value Theorem,  that $(|\theta_1|+ R^{-1/2})|\alpha_1|\ge \delta$, so that  
    $$|\alpha_1|\ge \min\(R^{-1/2}, \frac{\delta}{\e R+\sqrt R}\),$$
    where 
    \begin{equation}\label{eps}\e = \max(|\alpha_{21} - \pi/2|, |\beta_{12} + \pi/2|, |\alpha_{32} - \pi/2|, |\beta_{23} + \pi/2|).\end{equation}
    (Note that $\e = o(1)$ as $n\to +\infty$ and $|\theta_1| \le 2\e$.) 
    
    Therefore, since  $|q-a_1|=R| e^{i\alpha_1}-1|$,  and using Taylor's expansion of the sine function, 
    $$|q-a_1|\ge \frac12 \min\(\sqrt R,\frac{\delta}{\e + R^{-1/2}}\) .
    $$
    Similar estimates hold for $b_1$, and for $a_3$, $b_3$ using the function $\zeta_3$. From this we deduce that if 
    \begin{equation}
    \delta(t) = \min(|t-\zeta(a_1)|,|t-\zeta(b_1)|, |t-\zeta(a_3)|,|t-\zeta(b_3)|),\label{delta}
    \end{equation}
    then for any $q\in\ga_t\cap\partial B_R$ we have
    $$\min(|q-a_1|,|q-b_1|, |q-a_3|,|q-b_3|)\ge\frac12 \min\(\sqrt R, \frac{\delta(t)}{\e+R^{-1/2}}\).
    $$
    Therefore, using the particular boundary data we have on $\partial B_R$, if $q\in\ga_t\cap U_i$ then 
    $$ d(p_i,q) \le \exp{\(-c\min\(\sqrt R,\frac{\delta(t)}{\e + R^{-1/2}}\) \)},$$
    from which we deduce that
    \begin{equation}\label{lowslice}E(u,\ga_t) \ge d(p,p') - \exp{\(-c\min\(\sqrt R,\frac{\delta(t)}{\e + R^{-1/2}}\) \)},\end{equation}
    where $p=p_1$ and $p' = p_2$ in case a), where $p=p_2$ and $p'=p_3$ in case b), and where $p=p_1$ and $p'=p_3$ in case c). 

    If we integrate \eqref{lowslice} with respect to $t$,  and in view of  \eqref{Fred}, \eqref{delta},  we find
    $$E(u,B_R) \ge |a_1 - b_1| d_{12} +|a_3 - b_3| d_{23} - C\(\e + R^{-1/2} +Re^{-c\sqrt R}\).$$
    This is not enough to contradict \eqref{up2}, but we may now use \eqref{central} to improve the lower-bound. First, using again the fact that the transition angles tend to $\pm\pi/2$ as $n\to +\infty$,  we note that $u_1$ and $u_3$ tend to $(0,1)$ as $n\to +\infty$. This implies that, for $n$ large enough, 
    \begin{equation}
        |\mathcal{T}|\ge 2R_0 - \eta,\quad \text{where}\quad \mathcal{T} = \{t\in\R\mid \text{$\ga_t\cap\mathcal R_\lft\ne\emptyset$ and $\ga_t\cap\mathcal R_\rght\ne\emptyset$}\}.\label{T}
    \end{equation}
    For each $t\in \mathcal{T}$, $\ga_t$ is a broken line of length $2R + o(1)$, and from the definition of $\mathcal R_\lft$ and $\mathcal R_\rght$ there are points  $s_1$ $s_3\in\ga_t$  such that 
    \begin{equation}
    d(u(s_1),p_1) < \eta, \quad d(u(s_3),p_3)<\eta,\quad s_1,\ s_3\in B_{R_0}.\label{points}
    \end{equation}
    These two points divide $\ga_t$ in three portions $I_\lft$, $I_0$ and $I_\rght$. The energy of $u$ on $I_\lft$ and $I_\rght$ is bounded below using Lemma~\ref{lem20} while the energy on $I_0$ is bounded below using Lemma~\ref{lem21}. We find that
    $$ E(u,I_\lft) \ge d(p_1,s_1) - Ce^{-cR},\quad E(u,I_\rght) \ge d(s_3,p_3) - Ce^{-cR},\quad E(u,I_0) \ge d(s_1,s_3) + ce^{-CR_0}.$$
    Adding  we obtain
    \begin{equation} E(u,\ga_t) \ge d_{13} + ce^{-CR_0} - Ce^{-cR}.\label{lowfinal}\end{equation}
    We now integrate with respect to $t$ our lower bounds, using either \eqref{lowslice} or \eqref{lowfinal} according to whether $t\in \mathcal{T}$ or not. In view of \eqref{T} we find that 
    $$E(u,B_R) \ge R\(|e^{i\alpha_{21}}- e^{i\beta_{12}}| d_{12} + |e^{i\alpha_{32}}- e^{i\beta_{23}}| d_{23}\) +cR_0e^{-cR_0} - \Delta,$$
    where 
    $$\Delta = C\(\e + R^{-1/2} + Re^{-c\sqrt R} + Re^{-cR\eta}\).$$
    Since $\Delta\to 0$ as $n\to +\infty$, we obtain a contradiction with \eqref{up2} if $n$ is large enough, proving Theorem~\ref{main}.

\bibliographystyle{acm}
\bibliography{EP}

\end{document}